\documentstyle[11pt]{article}

\newfont{\msbm}{msbm10}

\def\U{{\cal U}}

\def\V{{\cal V}}

\def\L{{\cal L}}

\def\cinf{C^{\infty}}

\def\S{{\cal S}}

\def\N{\hbox{\msbm N}}

\def\Z{\hbox{\msbm Z}}

\def\R{\hbox{\msbm R}}

\def\C{\hbox{\msbm C}}

\def\H{{\cal H}}

\def\br{\bar{\R}}

\def\m{\mu}

\def\f{\phi}

\def\vf{\varphi}

\def\s{\sigma}

\def\h2{${\rm h}(2)$}

\newtheorem{lem}{Lemma}
\newtheorem{prop}{Proposition}
\newtheorem{cor}{Corollary}
\newtheorem{theo}{Theorem}
\newtheorem{deff}{Definition}

\def\bc{\begin{cor}}
\def\ec{\end{cor}}

\def\bt{\begin{theo}}
\def\et{\end{theo}}

\def\bd{\begin{deff}}
\def\ed{\end{deff}}

\def\bp{\begin{prop}}
\def\ep{\end{prop}}

\def\ba{\begin{eqnarray}}
\def\ea{\end{eqnarray}}

\def\be{\begin{equation}}
\def\ee{\end{equation}}

\def\reps{representations}

\newfont{\msbms}{msbm6} 

\def\ss{\S(\R^d)}

\def\B{{\cal B}}

\def\Ni{\hbox{\msbms N}}

\def\Ri{\hbox{\msbms R}}

\def\Bs{\hbox{\msbm B}}

\begin{document}

\title{Topics of Measure Theory on  Infinite Dimensional~Spaces}

\author{Jos\'e Velhinho}

\date{{Faculdade de Ci\^encias, Universidade da Beira 
Interior\\R. Marqu\^es D'\'Avila e Bolama,
6201-001 Covilh\~a, Portugal}\\{jvelhi@ubi.pt}}

\maketitle

\begin{abstract}
\noindent 
This short review is devoted to measures on infinite
dimensional spaces. We start by discussing product measures and projective techniques.
Special attention is paid to measures on linear spaces, and in particular 
to Gaussian measures.
Transformation properties of measures are considered, as well as 
fundamental results concerning the support of the measure.
\end{abstract}



\section{Introduction}
\label{Meddimeninf}

We present here a brief introduction to the subject of measures on infinite
dimensional spaces. The author's background is mathematical physics and quantum field theory,
and that is likely to be reflected in the text, 
but an  hopefully successful effort was made to produce a  review of interest to a broader audience. 
{We  have   references  \cite{Ya,GV,S} as our main inspiration.
Obviously, some important topics are not dealt with, and others are discussed
from a particular perspective, instead of another.
Notably, we do not discuss the perspective of abstract Wiener spaces, emerging from the works of
 Gross  and others \cite{Sa,Se,Pr,Gr}. Instead, we approach  measures in 
general linear spaces from the projective perspective (see below).}

For the sake of completeness, we include in Section \ref{geral} 
fundamental  notions and definitions from measure theory, with particular attention 
to the  issue of $\s$-additivity. We start by considering in Section~
\ref{espacoproduto} 
the infinite product of a  family of  probability measures.
In Section \ref{limproj} we consider  projective techniques, which play an
important role in  applications {(see  e.g., \cite{AL2,Le} for applications to gauge theories
and gravity)}.
Sections \ref{esplin} to \ref{mgauss} are devoted to measures on
 infinite dimensional  linear spaces.
In Section \ref{B-M-I} results concerning the support of the measure are presented,
which partly justify, in this context, the interest of nuclear spaces and their 
(topological) duals. The 
 particular case of 
Gaussian measures is considered in 
Section \ref{mgauss}. 
{There are  of course several possible approaches to the issue of  measures in 
infinite dimensional linear spaces, and to Gaussian measures in particular, including the well
known and widely used framework of Abstract Wiener Spaces or other approaches
working directly with Banach spaces (see, e.g., \cite{Sat,K,Boga}).
We follow here the approach of Ref. \cite{Ya}, taking advantage of the facts
that the algebraic dual of any linear space is a projective limit (of finite dimensional spaces)
and that any consistent family of measures defines a measure on the projective limit.}
In Section \ref{srepci} we present the main  definitions and some fundamental results 
concerning  transformation properties of measures, discussing briefly
quasi-invariance and ergodicity.
Finally, in~Section \ref{temp} we consider in particular measures on the space
of tempered distributions.

Generally speaking, and except when explicitly stated otherwise,
we consider  only finite (normalized) measures.
(A notable exception is the  Lebesgue measure on $\R^n$, to which 
we refer~occasionally.)


\section{Measure Space}
\label{geral}
We review in this  section some fundamental aspects of measure theory,
focusing (although not 
exclusively) on finite measures. A very good presentation  of these
subjects can be found {in \cite{RS1,KF,B1,B,Cohn}.}

\bd
\label{mdef1} Given a set $M$, a family $\cal F$ of subsets
of $M$ is said to be a  (finite) algebra if it is closed under the operations of
taking the complement and  finite unions, i.e., if $B\in{\cal F}$ implies
$B^{\rm c}\in{\cal F}$ and $B_1\in{\cal F},\ldots ,B_n\in{\cal F}$ implies
$\cup _iB_i\in {\cal F}$. 
\ed
\bd
\label{mdef2} A  non-negative real function $\mu$ on an algebra
$\cal F$ is said to be a {measure} if for any finite set
 of mutually disjoint elements $B_1,B_2,\ldots ,B_n$ of 
 ${\cal F}$ 
{\rm (}$B_i\cap B_j=\emptyset$ for $i\neq j$\/{\rm )} 
the {following additivity} condition is~satisfied:
\be
\mu \bigl({\cup}_iB_i\bigr)=\sum _i\mu \,(B_i) .
\ee
\ed
{Particularly important is the notion of  measures on   $\s$-algebras,
in which case the measure is required to satisfy the so-called $\sigma $-additivity
condition.}
\bd
\label{mdef3} 
Given a set $M$, a family $\cal B$ of subsets of
 $M$ is said to be a  $\sigma $-algebra if it is closed under 
complements and countable unions, i.e.,  $B\in{\cal B}$ implies
$B^{\rm c}\in{\cal B}$ and $B_i\in{\cal B}$, $i\in\N$, implies
$\cup_i^{\infty}B_i\in {\cal B}$. 
The pair 
$(M,{\cal B\,})$ is called a  measurable space and the  elements of $\cal B$ 
are called measurable sets.
\ed
It is obvious that for any   measurable space
$(M,{\cal B\,})$ the $\sigma $-algebra $\B$ contains $M$ and the empty set,
and it is also closed under countable intersections.
Another  operation of 
interest in a  $\s$-algebra (or~finite algebra) is the symmetric difference
of  sets $A\triangle B:=(A\setminus B)\cup(B\setminus A)=
(A\cup B)\setminus(A\cap B)$.
\bd
\label{mdef4} 
{Given a measurable space 
$(M,{\cal B\,})$, a function $\mu:\B\rightarrow[0,\infty]$, with 
$\mu(\emptyset)=0$, is said 
to be a  
measure if it satisfies the
$\sigma $-additivity property, i.e., if for any sequence of mutually disjoint
 measurable sets  $\{B_i\}_{i \in \Ni}$ 
one has
\be
\mu \bigl({\cup}_i^{\infty}B_i\bigr)=\sum _i^{\infty}\mu \,(B_i) ,
\ee
where the right hand side denotes either the sum of the series or infinity, 
in case the sum does not converge.
The~structure $(M,{\cal B},\mu )$
is called a measure space. The measure is said to be finite if $\m(M)<\infty$,
and normalized if $\m(M)=1$, in which case $(M,{\cal B},\mu )$
is  said to be a   probability space.}
\ed
An  important property following from $\sigma $-additivity is the following.
\bt
\label{mteo1}
Let $\mu $ be a   $\s$-additive finite measure and
$B_1\supset B_2\supset \ldots$
a decreasing sequence of measurable sets.~Then
\be
\mu \bigl(\cap_nB_n\bigr)=\lim _{n \to \infty }\mu \,
(B_n).
\ee
also,
\be
\mu \bigl(\cup_nA_n\bigr)=\lim _{n \to \infty }\mu \,
(A_n),
\ee
for any increasing sequence  $A_1\subset A_2 \subset \ldots $ 
of measurable sets.
\et
Let us consider the problem of  the extension  of measures on finite
algebras to  ($\sigma $-additive) measures on $\sigma $-algebras.
Note first that given any  
family $\cal A$ of subsets of a set $M$ there is a 
 minimal  $\sigma$-algebra
containing $\cal A$. We will denote  this $\sigma$-algebra by
$\Bs({\cal A})$, the $\sigma$-algebra generated by $\cal A$.

\bt[{Hopf \cite{hopf}}]
\label{mteo3} A finite measure $\mu $ on an algebra $\cal F$ can be 
extended to  $\sigma $-additive finite measure on the $\sigma $-algebra 
$\Bs({\cal F})$ if and only if for any given 
decreasing sequence 
$B_1\supset B_2\supset \ldots $ of elements of $\cal F$ the condition
$\lim _{n \to \infty}\mu \,(B_n)>0$ implies $\cap_nB_n\neq \emptyset$.
\et
\bt
\label{mteo2} If it exists, the extension of a finite measure on $\cal F$ 
to a $\sigma $-additive finite measure on $\Bs({\cal F})$  is unique.
\et
Among non-finite measures,  
so-called $\sigma$-finite measures are particularly important.
\bd
\label{sfinite} A measure is said to be $\sigma$-finite if the measure
space $M$ is a countable  union of mutually disjoint
 measurable sets, 
each of which with  finite measure. 
\ed
The Lebesgue measure on $\R ^n$ is of course $\s$-additive and $\sigma$-finite.
\bd
\label{Borel}
Let $(M,\tau)$ be a  topological space,  $\tau$ being the family of open sets.
The $\s$-algebra $\Bs(\tau)$ generated by open sets is called a
 Borel $\s$-algebra. The measurable space $\bigl(M,\Bs(\tau)\bigr)$
is said to be Borel (with respect to $\tau$). A  
{measure} on 
$\bigl(M,\Bs(\tau)\bigr)$ is called a Borel measure.
\ed
Except when explicitly said otherwise,  $\R^n$ and 
$\C^n$ are considered to be equipped with the  usual topology and
corresponding  Borel $\sigma $-algebra.  
\bd
\label{regularBorel}
A  Borel measure $\m$ is said to be  regular if for any 
 Borel set $B$ one has:
\ba
\label{regularBorel1}
\m(B)&=&\hbox{\rm inf}\, \{\m(O)\ |\ \ O\supset B,\ O\ 
open\}\nonumber \\
\label{regularBorel2}
&=&\hbox{\rm sup}\, \{\m(K)\ |\ \ K\subset B,\ K\ 
compact\ and\ Borel\}.
\ea
\ed
\bp
\label{geral2}
Any Borel measure on a separable and complete metric space is
regular.
\ep
\bd
\label{mdef8}
Let $(M,\B,\m)$ be a measure space and  ${\cal N}_{\m}:=\{B\in\B\ |\ \ 
\m(B)=0\}$ the  family  of  zero measure sets. Two sets
$B_1,B_2\in\B$ are said to be equivalent modulo zero measure sets, 
$B_1\sim B_2$, if and only if $B_1\,\triangle \,B_2\in{\cal N}_{\m}$.
\ed
The family ${\cal N}_{\m}$ of zero measure sets is an ideal
on the ring of  measurable sets  $\B$ defined by the  operations
$\triangle$ and $\cap$, and therefore the quotient $\B/{\cal N}_{\m}$ 
is also a ring.
It is straightforward  to check that the measure is well   defined 
on $\B/{\cal N}_{\m}$. From the strict measure theoretic  point of view, 
the fundamental objects are the
elements of $\B/{\cal N}_{\m}$, and   {naturally defined} transformations between 
measure spaces $(M,\B,\m)$ and $(M',\B',\m')$ are maps between
the quotients $\B/{\cal N}_{\m}$ and $\B'/{\cal N}_{\m'}$.
\bd
\label{mdef9}
Two measure spaces $(M,\B,\m)$ and $(M',\B',\m')$ are said to be isomorphic
if there exists a bijective transformation between 
$\B/{\cal N}_{\m}$ and $\B'/{\cal N}_{\m'}$, mapping $\m$ into $\m'$.
\ed
In the above sense,  zero measure sets are  irrelevant. (When the measure is defined in a topological space, a more restricted notion of {\em support of the measure} is sometimes adopted, namely the  smallest closed set with full measure. We do not adhere to that definition of support.)
\bd
\label{defbm3}
Let $(M,\B,\mu )$ be a measure space. A (not necessarily measurable) 
subset $S\subset M$ is said to be a support of the measure 
if any measurable subset in its complement has zero measure, i.e.,
$Y\subset S^{\rm c}$ and $Y\in\B$ implies $\m(Y)=0$.
\ed
Given a measurable space $(M,\B)$ and a subset $N\subset M$,
let us consider the $\s$-algebra of measurable subsets of  $N$,
\be
\label{algebra subconjunto}
\B\cap N:=\{B\cap N\ |\ \ B\in\B\}.
\ee

If $(M,\B,\m)$ is a  measure space 
 and  
$N\subset M$ is  measurable, there is a naturally 
defined   {measure} $\m_{_{|N}}$ on $(N,\B\cap N)$, by 
restriction of $\m$ to $\B\cap N$, $\m_{_{|N}}(B\cap N):=\m(B\cap N)$, 
$\forall B\in\B$. The restriction of the measure is also well defined 
for subsets $S\subset M$  supporting the measure, even if 
 $S$ is not a measurable set. In this case we have
$\m_{_{|S}}(B\cap S)=\m(B)$. One can in fact show the following \cite{Ya}.
\bp
\label{suportegeral}
If $\m$ is a  {measure} on $(M,\B)$ and $S$ is a support of the measure
then  $\m_{_{|S}}(B\cap S):=\m(B)$, $\forall B\in\B$, defines a
{($\s$-additive)} measure on $(S,\B\cap S)$. The measure spaces $(M,\B,\m)$ and 
$(S,\B\cap S,\m_{_{|S}})$ are~isomorphic.
\ep
[The measure on $S$ is well defined, since 
$B_1\cap S=B_2\cap S$ implies $(B_1\,\triangle\,B_2)\cap S=\emptyset$, which
in turns leads to  $\m(B_1\,\triangle\, B_2)=0$, given that $S$  supports the measure.]
\bd
\label{mdef6} A transformation  $\vf:M_1\to M_2$ between two  measurable spaces
 $(M_1,{\cal B}_1)$ and $(M_2,{\cal B}_2)$ is said to be measurable
if $\vf^{-1}\B_2\subset\B_1$, i.e., if
$\vf^{-1}B\in{\cal B}_1$, $\forall B\in {\cal B}_2$, where $\vf^{-1}B$ is
the preimage of $B$.
\ed
Given a measurable transformation $\vf :M_1\rightarrow M_2$ between
measurable spaces $(M_1,{\cal B}_1)$ and $(M_2,{\cal B}_2)$, one gets a map 
$\tilde \vf :{\cal B}_2\rightarrow{\cal B}_1$, defined by 
$\tilde\vf(B)=\vf^{-1}B$. If $\mu $ is a   {measure} on
$(M_1,{\cal B}_1)$, the composition map  $\mu\circ\tilde\vf$ is therefore
 a  
{measure} on $(M_2,{\cal B}_2)$, defined by:
\be
\label{II1.5.3}
(\mu\circ\tilde\vf)(B)=\mu (\varphi ^{-1}B),\ \forall B\in {\cal B}_2.
\ee

This measure is usually called the push-forward of $\mu $ with respect to $\varphi $. 
[Given that a measure
$\m$ on $(M,\B)$ is in fact a function on $\B$, the measure $\m\circ\tilde\vf$ is
actually the pull-back of  $\m$ by $\tilde\vf$; we will use however the
usual expression ``push-forward''.]
Besides $\m\circ\tilde\vf$, we will use also the  alternative notations
$\vf_*\,\mu$ and  $\mu_{\vf}$ to denote the push-forward of a measure.

Measure  theory  is naturally connected to  integration. 
From this  point of view, the (in~general~
complex) measurable functions $f:M\to\C$ are particularly 
important, in a measure space  $(M,\B,\m)$. 
More precisely, the  relevant objects are
 equivalence classes  of measurable functions.
\bd
\label{mdef10}
Given a measure space $(M,\B,\m)$, a condition $C(x)$, 
$x\in M$, is said to be satisfied almost everywhere if the set: 
$$\{x\in M\ |\  C(x)\ {is\ false}\}$$ 
is contained in a  zero measure set.
\ed
\bd
\label{mdef11}
Two   measurable complex functions $f$ and $g$ on a measure space are said to
be  equivalent
if the  condition $f(x)=g(x)$ is satisfied almost everywhere.
\ed
The set of  equivalence  classes of measurable functions 
is naturally a linear space. 
With a finite measure $\m$,  the  {integral} 
defines a family of norms, by: 
\be
\|f\|_p:=\left(\int |f|^pd\m\right)^{1/p},
\ee 
with $p\geq 1$.
With the norm $\|\ \|_p$, the linear space of equivalence  classes of measurable functions  
is denoted by $L^p(M,\m)$. The space  $L^p$ is defined analogously  for non-finite measures,
considering only functions such that  the integral over the whole space is finite, 
$\int |f|^p d\m<\infty$. Let us recall still that in the  particular  case $p=2$ the norm
comes from an inner  product, $(f,g)=\int f^*gd\m$, and therefore the space 
$L^2(M,\m)$   of  (classes of) 
 square integrable (complex) functions on  $(M,\B,\m)$ is an
inner product space.
In this context, the  interest of
$\s$-additive measures is rooted in the crucial fact that the $L^p$ 
spaces associated with these measures are complete.
{Except when explicitly stated (namely when the question
of $\s$-additivity is explicitly concerned),
we will drop the qualifier ``$\s$-additive'' when referring to measures
on $\s$-algebras.}

The next result, which follows from the definition of 
{integral},
generalizes the usual change of~variables.
\bp
\label{mudancavariavel}
Let $(M,{\cal B},\mu )$ be a measure space, $(M',\B')$ a
measurable space and $\varphi :M \to M'$ a measurable transformation.
Consider the measure space $(M',{\cal B}',\mu _{\varphi })$,
where $\mu _{\varphi }$ denotes the  push-forward with respect to $\varphi $.
Then, for any $\mu _{\varphi }$-integrable function $f:M'\rightarrow \C$, the 
function  $f\circ \varphi :M\rightarrow \C$ is 
integrable with respect to~$\mu $~and:
\be
\label{mudvar}
\int _M(f\circ \varphi )d\mu =\int _{M'}fd\mu _{\varphi }.
\ee
\ep

\section{Product Measures}
\label{espacoproduto}
Let $\bigl\{(M^1,\B^1,\m^1),\ldots,(M^n,\B^n,\m^n)\bigr\}$ be a finite set of
probability spaces. Consider the Cartesian~product:
\be
\label{prod1}
M_{n}:=\prod_{k=1}^n M^{k},
\ee
the projections:
\be
\label{prod2}
p_{n}^{k}:M_{n}\to M^{k}
\ee
and the $\sigma$-algebra of subsets of $M_{n}$:
\be
\label{prod3}
\B_{n}:=\Bs\biggl(\bigcup_{k=1}^n
(p_{n}^{k})^{-1}\B^{k}\biggr) .
\ee

The measurable product space  of the spaces 
$\bigl\{(M^1,\B^1),\ldots,(M^n,\B^n)\bigr\}$ is the pair
$(M_{n},\B_{n})$. Note that $\B_{n}$ is the smallest
$\sigma$-algebra such that all projections
$p_{n}^{k}$ are measurable.

The $\sigma$-algebra $\B_{n}$  obviously contains the  Cartesian products of measurable sets 
 $\omega^k\in\B^k$, $k~=~1,\ldots,n$, i.e., $\B_{n}$
contains all sets of the form:
\be
\label{prod4}
(\omega^1,\ldots,\omega^n)=:\prod_{k=1}^n\omega^k,\ \
\omega^k\in\B^k,\ k=1,\ldots,n.
\ee

It is a classic result that there exists a unique probability measure 
 $\m_n$ in $(M_n,\B_n)$ such that:
\be
\label{prod5}
\m_n\biggl(\prod_{k=1}^n\omega^k\biggr)=\prod_{k=1}^n\m^k(\omega^k) ,
\ee
which is called the  product measure and is represented by:
\be
\label{prod6}
\m_n=\prod_{k=1}^n\m^k.
\ee

Let us consider now the  infinite product, not necessarily countable. 
As we will see immediately, the existence and  uniqueness of  the product measure
continue to take place.
\bd
\label{defprod1}
Let $\bigl\{(M^{\lambda},\B^{\lambda})\bigr\}_{\lambda\in\Lambda}$ be a  
family of measurable spaces labeled by a set $\Lambda$ and let
$M_{\Lambda}$ be the Cartesian product of the spaces $M^{\lambda},\ \lambda\in
\Lambda$. For each $\lambda\in\Lambda$ let $p_{\Lambda}^{\lambda}$ be the 
projection from $M_{\Lambda}$ to $M^{\lambda}$.  The~measurable product space of the  family $\bigl\{(M^{\lambda},\B^{\lambda})\bigr\}_
{\lambda\in\Lambda}$ 
is defined as the pair $(M_{\Lambda},\B_{\Lambda})$, where:
\be
\label{prod11}
\B_{\Lambda}:=\Bs\biggl(\bigcup_{\lambda\in\Lambda}
(p_{\Lambda}^{\lambda})^{-1}\B^{\lambda}\biggr)
\ee
is the smallest $\sigma$-algebra such that all projections
$p_{\Lambda}^{\lambda}$ are measurable.
\ed
Consider now a family 
$\bigl\{(M^{\lambda},\B^{\lambda},\mu^{\lambda})\bigr\}_{\lambda\in\Lambda}$ 
of probability spaces and let  $\cal L$ be the family
of finite subsets of  $\Lambda$. For each $L\in{\cal L}$ 
let us consider the (finite) product   probability space  $(M_L,\B_L,\m_L)$ 
defined as above, i.e.,
\be
\label{prod12}
M_{L}=\prod_{\lambda\in L}M^{\lambda},
\ee
\be
\label{prod13}
\B_{L}=\Bs\biggl(\bigcup_{\lambda\in L}
(p_{L}^{\lambda})^{-1}\B^{\lambda}\biggr)
\ee
(where $p_{L}^{\lambda}$ is the projection from $M_L$ to $M^{\lambda}$) and,
\be
\label{prod14}
\m_L=\prod_{\lambda\in L}\m^{\lambda}.
\ee

Consider still the  natural measurable  projections,
\be
\label{prod15}
p_{L,\Lambda}:M_{\Lambda}\to M_L .
\ee

The following result can be found in  \cite{Ya}.
\bt
\label{teoprod1}
There is a unique  
($\s$-additive) probability measure
$\m_{\Lambda}$ in $(M_{\Lambda},\B_{\Lambda})$ such that:
\be
\label{prod16}
(p_{L,\Lambda})_*\,\m_{\Lambda}=\m_L,\ \forall L\in{\cal L}.
\ee
\et
The measure  defined by this  theorem is called the   product measure.

\vspace{12pt}
\noindent {\bf Example:} A simple but important example of a  product measure on  
 an infinite dimensional space is the following, which generalizes the notion of product Gaussian measures
in $\R^n$. Consider the countable
family of measurable spaces 
$\bigl\{(M^k,\B^k)\bigr\}_{k\in\Ni}$, where, for each $k$, $(M^k,\B^k)$ coincides with 
$\R$ equipped with the Borel
$\s$-algebra. The  measurable product space is the space 
$\R^{\Ni}$ of all real sequences:
\be
\label{prodex1}
x=(x_k)=(x_1,x_2,\ldots),
\ee
equipped with the smallest $\s$-algebra such that all projections $x\mapsto x_k$
are measurable. Let us consider in each of the spaces $\R$ of the
family the same Gaussian measure of covariance $\rho\in\R^+$, i.e.,
\be
\label{prodex2}
d\m^k(x_k)=e^{-x_k^2/2\rho}\frac{dx_k}{\sqrt{2\pi\rho}},\ \ \forall k\in\N.
\ee

According to the Theorem \ref{teoprod1}, the product measure, here denoted by
$\m_{\rho}$,
\be
\label{prodex3}
d\m_{\rho}(x)=\prod_{k=1}^{\infty}e^{-x_k^2/2\rho}\frac{dx_k}{\sqrt{2\pi
\rho}},
\ee
is uniquely determined  by its value  on the sets of the form
$\prod_{k\in\Ni}\omega^k$, 
where only for a finite subset of  $\N$ the Borel sets (in $\R$)
$\omega^k$ differ from $\R$.

\smallskip
Obviously, the above example can be generalized  for any  infinite sequence
of probability measures on $\R$, not necessarily identical. 
The correspondent of the  Lebesgue measure, ``$\prod_{k=1}^{\infty}dx_k$'', 
however, does not exist, i.e., the infinite product of Lebesgue measures
in $\R$ does not define a measure.

Given any  product measure, defined by a not necessarily countable family of probability spaces, 
it is also trivial to determine the measure of sets of the form:
\be
\label{prod18}
Z\bigl(\{\omega^{\lambda}\}\bigr):=\prod_{\lambda\in\Lambda}\omega^{\lambda},
\ \ \omega^{\lambda}\in\B^{\lambda},
\ee
where only for a countable subset  $\Lambda_0=\{\lambda_i\}_{i\in\Ni}
\subset\Lambda$ the sets $\omega^{\lambda}$ differ from $M^{\lambda}$. 
Since it is a typical argument in measure theory, we present it next in some detail.
Let us start by showing that the sets~(\ref{prod18}) are measurable. 
Consider the  finite subsets of $\Lambda_0$,  $L_n:=\{\lambda_1,\ldots,
\lambda_n\},\ n\in\N$, and let $Z_n$ be the sets  defined as in 
(\ref{prod18}), but where $\omega^{\lambda_k}$ ($\lambda_k\in\Lambda_0$) is replaced by
 $M^{\lambda_k}$ for $k>n$. It is clear that:
 \vspace{-6pt}
\be
\label{prod19}
Z_n=p_{L_n,\Lambda}^{-1}\prod_{k=1}^n\omega^{\lambda_k},
\ee
and it follows that $Z_n$ is measurable. Since $Z_{n'}\subset Z_n$
for $n'>n$, the sets $\{Z_n\}_{n\in\Ni}$ form a decreasing
sequence of measurable sets. 
The intersection $\cap_{\Ni}Z_n$
is therefore measurable, since $\B_{\Lambda}$ is a  $\s$-algebra. 
But $Z\bigl(\{\omega^{\lambda}\}\bigr)$ coincides precisely with  
$\cap_{\Ni}Z_n$. 
Invoking the $\s$-additivity of the measure we then get from theorem 
\ref{mteo1} and (\ref{prod16}):
 \vspace{-6pt}
\ba
\label{prod20}
\m_{\Lambda}\Bigl(Z\bigl(\{\omega^{\lambda}\}\bigr)\Bigr)&=&
{\lim}_{n\to\infty}\m_{\Lambda}(Z_n) \nonumber \\
\label{prod21}
&=&{\lim}_{n\to\infty}\m_{\Lambda}\biggl(p_{L_n,\Lambda}^{-1}
\prod_{k=1}^n\omega^{\lambda_k}\biggr) \nonumber \\
\label{prod22}
&=&{\lim}_{n\to\infty}\m_{L_n}\biggl(\prod_{k=1}^n\omega^{\lambda_k}\biggr) 
\nonumber \\
\label{prod23}
&=&{\lim}_{n\to\infty}\prod_{k=1}^n\m^{\lambda_k}(\omega^{\lambda_k}).
\ea

\section{Projective Limits}
\label{limproj}

We present in this section the notion of measurable  projective limit space.

Let us start by recalling that a set
 $\cal L$ is said to be  partially ordered if it is equipped with a partial order 
relation, i.e., there is a binary relation
 ``$\geq $'' such that:
\begin{enumerate}
\item[(1)]  (reflexivity) $L\geq L,\ \forall L \in {\cal L}$
\item[(2)]  (transitivity) $L\geq L'\ {\rm and}\ L'\geq L''
\Rightarrow L\geq L''$
\item[(3)]  (anti-symmetry) $L\geq L'$ {\rm and} $L'\geq L
\Rightarrow L=L' .$
\end{enumerate}

Recall still that a set $\cal L$, partially ordered with respect to the partial order 
relation ``$\geq $'', is said to be directed if
$\forall L',L'' \in
{\cal L}$ there exists $L \in {\cal L}$ such that $L\geq L'$ and $L\geq L''$.
\bd
\label{deflimproj1}
Let $\cal L$ be a directed set and $\{M_L\}_{L\in {\cal L}}$ a family of sets labeled by
 $\cal L$. Suppose that for each pair
 $L',L$ such that $L'\geq L$ there are  surjective maps:
\be
\label{lim1}
p_{L,L'}:\ M_{L'}\rightarrow M_L
\ee
satisfying:
\be
\label{lim2}
p_{L,L'}\circ p_{L',L''}=p_{L,L''},\ {\rm for}\ \,L''\geq L'\geq L.
\ee

The family $\{M_L,p_{L,L'}\}_{L,L'\in {\cal L}}$ of sets $M_L$ and
maps $p_{L,L'}$ is called a  projective family. 
\ed
>From now on, the maps \ $p_{L,L'}$ of the projective family will be called projections. Let us consider the Cartesian product of the sets  $M_L$:
\be
\label{lim3}
M_{\cal L}:=\prod _{L\in {\cal L}}M_L,
\ee
and denote its generic element by $(x_L)_{L\in {\cal L}},\ 
x_L \in M_L$. 
\bd
\label{deflimproj2}
The  projective limit of  the family $\{M_L,p_{L,L'}\}_{L,L'\in {\cal L}}$ 
is the subset  $M_{\infty }$ of the  Cartesian product $M_{\cal L}$ 
defined by:
\be
\label{lim4}
M_{\infty }:=\Bigl\{ (x_L)_{L\in {\cal L}}\in M_{\cal L}\ |\ \ 
L'\geq L \Rightarrow p_{L,L'}\,x_{L'}=x_L \Bigr\}.
\ee
\ed
The  projective limit  is therefore formed by  consistent families of
 elements $x_L \in M_L$, in the sense that $x_L$ is  defined by
$x_{L'}$, for $L'\geq L$.
\bd
\label{deflimproj3}
A family $\bigl\{ (M_L,\B _L),p_{L,L'}\bigr\}_{L,L'\in{\cal L}}$
is said to be a measurable projective family if each pair 
$(\,M_L,\B _L\,)$ is a measurable space and if $\{M_L,p_{L,L'}\}$ is
a  projective family such that all projections $p_{L,L'}$ are~measurable.
\ed
Given a measurable projective family, the structure of measurable space in the
projective limit $M_{\infty}$ is defined as follows.
Let $\B _{\cal L}$ be the  product $\sigma $-algebra defined in the previous section, 
i.e., $\B _{\cal L}$ is the smallest $\sigma $-algebra of subsets of the product space
$M_{\cal L}$ such that all the projections from
 $M_{\cal L}$ to $M_L$ are measurable (note that, with respect to the product, the spaces
 $M_L$ play here the role of the spaces 
$M^{\lambda }$ of the previous section). Let us consider the
$\s$-algebra $\B _{\infty}$ of subsets of $M_{\infty}$ given by:
\be
\label{lim5}
\B _{\infty}:=\B _{\cal L}\cap M_{\infty}=\{B\cap M_{\infty}\ |\ \ 
B\in \B _{\cal L}\}.
\ee

The family $\B _{\infty}$ is closed under countable unions, since:
\be
\textstyle{\bigcup _n}\Bigl(B_n\bigcap M_{\infty}\Bigr)=\Bigl(\bigcup _nB_n
\Bigr)\bigcap M_{\infty}.
\ee

Let us also show that $\B _{\infty}$ is closed under the operation of taking the complement, 
i.e., that~$M_{\infty}\backslash (B\cap M_{\infty})\in \B _{\infty},
\ \forall B\in \B _{\cal L}$. Taking $M_{\infty}$ and $B\cap M_{\infty}$
as subsets of $M_{\cal L}$ we get:
\ba
\label{lim7}
M_{\infty}\backslash (B\cap M_{\infty}) & = & 
M_{\infty}\cap\bigl(M_{\cal L}\backslash (B\cap M_{\infty})\bigr)\nonumber \\
\label{lim8}
& = & M_{\infty}\cap\bigl((M_{\cal L}\backslash B)\cup (M_{\cal L}
\backslash M_{\infty})\bigr)\nonumber \\
\label{lim9}
& = & M_{\infty}\cap\ (M_{\cal L}\backslash B) ,\nonumber
\ea
which proves the statement, since $M_{\cal L}\backslash B\in \B_{\cal L}$.
It follows that $\B _{\infty}$ as defined above is  indeed a $\sigma $-algebra.
\bd
\label{deflimproj4}
The pair $(M_{\infty},\B _{\infty})$ is called the  measurable projective limit of the 
measurable projective family 
$\bigl\{ (M_L,\B _L),p_{L,L'}\bigr\}_{L,L'\in{\cal L}}$.
\ed
Let $\pi_L$ be the  projection from $M_{\cal L}$ to $M_L$ and $p_L$ the restriction
of $\pi_L$ to $M_{\infty}$, i.e.,
\be
\label{lim10}
p_L=\pi _L\circ i_{\infty},
\ee
where $i_{\infty}$ is the inclusion of $M_{\infty}$ in $M_{\cal L}$. 
Since the maps {$\pi_L$ and $i_{\infty}$} are measurable by construction,
 $p_L$~is measurable $\forall L\in {\cal L}$. The 
 consistency conditions that define  $M_{\infty}$ are equivalent to:
\be
\label{lim11}
p_L=p_{L,L'}\circ p_{L'}\, ,\ \ \forall L,L':L'\geq L,
\ee
which in particular shows that:
\be
\label{lim11a}
{\cal F}_{\infty}:=\bigcup_{L\in{\cal L}}p_L^{-1}\B_L
\ee
is an algebra. The algebra ${\cal F}_{\infty}$ is formed by all the sets of the type
$p_L^{-1}B_L$, $B_L\in\B_L$, $L\in{\cal L}$,
which are called cylindrical sets.
One can further show that:
\be
\label{lim11b}
\B_{\infty}=\Bs({\cal F}_{\infty}) ,
\ee
and it follows that $\B_{\infty}$ is the smallest $\s$-algebra such that all
projections $p_L:M_{\infty}\rightarrow M_L$ are measurable~\cite{Ya}.

Suppose now that we are given a measure $\mu $ on $(M_{\infty},\B _{\infty})$.
The push-forward:
\be
\label{lim12}
\mu _L:=
(p_L)_*\,\mu
\ee
of $\mu$ by $p_L$ is a measure on $(M_L,\B _L)$. Explicitly:
\be
\label{lim12a}
\mu _L(B_L)=\mu \bigl(p_L^{-1}B_L\bigr),\ \ \forall B_L\in \B_L.
\ee

From (\ref{lim11}) it follows that the family of measures  
$\{\mu _L\}_{L\in {\cal L}}$ satisfy the self-consistency conditions:
\be
\label{lim13}
\mu _L=
(p_{L,L'})_*\,\mu_{L'}\,,\ \ \forall L,L':L'\geq L.
\ee

The problem of introducing a measure on a projective limit space is the inverse problem, 
i.e., we~look to define a measure on
$(M_{\infty},\B _{\infty})$ starting from a self-consistent family
of measures $\{\mu _L\}$. 

Note that given a self-consistent family
$\{\mu _L\}$ one can always define, by means of (\ref{lim12a}), an additive
measure $\mu$, called cylindrical, in ${\cal F}_{\infty}$. 
So, the problem consists  in the extension  of additive measures on
${\cal F}_{\infty}$ to  $\s$-additive measures on $\Bs({\cal F}_{\infty})$. 
An  important case where the  cylindrical measure can be extended to a 
$\sigma$-additive measure is that of the product measure of probability measures,  discussed in the previous section.  In fact, the product space can be seen as the  
projective limit of the family of finite products. In general, the existence
of measure on $(M_{\infty},\B_{\infty})$ depends  on 
topological conditions on the  projective family. 
Another particularly interesting situation where the  extension is ensured is the following
 \cite{Ya,AL2}.
\bd
\label{defmeba1}
A projective family $\{M_L,p_{L,L'}\}_{L,L'\in \L}$ of
 compact Hausdorff spaces is said to be a  compact  Hausdorff family if all the  projections
$p_{L,L'}$ are continuous.
\ed
One can show that the
projective limit of a compact  Hausdorff family is a 
compact  Hausdorff space, with respect to the  topology induced from the 
 Tychonov topology (in the  product space (\ref{lim3}))~\cite{MM}. 
\bt
\label{teomeba2}
Let $\{M_L,p_{L,L'}\}_{L,L'\in \L}$ be a   compact Hausdorff projective family.
Any self-consistent family of 
regular Borel probability measures  $\{\m_L\}_{L\in\L}$ in the family of spaces 
$\{M_L\}_{L\in\L}$ defines a  regular  Borel probability measure on 
the projective  limit $M_{\infty}$.
\et
Let us conclude this section with the notion of cylindrical functions
and a typical application of $\sigma$-additivity,
analogous to the result (\ref{prod23}) of the previous section.
Let us suppose then that we are given a 
{measure} $\mu$ on 
$(M_{\infty},\B_{\infty})$
and let  $\{\mu_L\}$ be the corresponding self-consistent family of measures
in the spaces $M_L$. Given an integrable function
$F$ on $M_{L_0}$, one gets by  pull-back an integrable function
$F\circ p_{L_0}$ on
$M_{\infty}$. Functions of this type are called cylindrical and they are 
the  simplest integrable functions on $M_{\infty}$. 
 From Proposition \ref{mudancavariavel} we get:
\be
\label{lim19}
\int _{M_{\infty}} \bigl(F\circ p_{L_0}\bigr)d\mu =\int _{M_{L_0}}Fd\m_{L_0}.
\ee

As a typical  example of the   construction of a non-cylindrical measurable set 
whose measure is trivially   determined, let us consider a countable subset
 ${\cal L}_0$ of ${\cal L}$, i.e., 
${\cal L}_0=\{L_1,L_2,\ldots\}$, with $L_{n+1}\geq L_n$ and let
$\{B_n\in \B_{L_n}\}_{n\in\Ni}$ be a sequence such that:
\be
\label{lim16}
p_{L_n,L_{n+1}}^{-1}B_n\subset B_{n+1}.
\ee
It is then clear from (\ref{lim11}) that:
\be
\label{lim17}
p_{L_n}^{-1}B_n \subset p_{L_{n+1}}^{-1}B_{n+1},
\ee
and $\bigl\{p_{L_n}^{-1}B_n\bigr\}$ is therefore an increasing sequence of 
cylindrical sets. 
The  union  of the sets $p_{L_n}^{-1}B_n$ is a measurable set which is in general 
non-cylindrical (it may be cylindrical if all sets
 $p_{L_n}^{-1}B_n$ coincide after some order).
 From Theorem \ref{mteo1} one therefore gets:
\be
\label{lim18}
\mu \Bigl({\textstyle \bigcup _n}p_{L_n}^{-1}B_n\Bigr)=\lim _{n\to \infty}
\mu _{L_n}\bigl(B_n\bigr).
\ee

\smallskip
\noindent {\bf Example:} The space 
known as the Bohr compactification of the line 
admits a projective characterization as follows (see \cite{bohr,SigB} for details).
For arbitrary $n\in\N$, let us consider  sets $\gamma=\{k_1,\ldots,k_n\}$ 
 of   real numbers $k_1,\ldots,k_n$,   such that the condition:
\be
\sum_{i=1}^n m_i k_i=0, \ \ \ m_i\in\Z,
\ee
can only be satisfied with $m_i=0$, $\forall i$. 
(These are of course sets  of linearly independent real numbers, with respect to the field of rationals.)
For each set $\gamma=\{k_1,\ldots,k_n\}$, let us consider the 
subgroup of $\R$ freely generated by $\gamma$:
\be
G_{\gamma}:=\left\{\sum_{i=1}^n m_i k_i,\ m_i\in\Z\right\}.
\ee

Let now  $T$ denote the group of unitaries in the complex plane, and
for each  $\gamma$ consider the  group $\R_{\gamma}$ of all group morphisms from 
$G_{\gamma}$ to $T$,
\be
\R_{\gamma}:={\rm Hom}[G_{\gamma},T].
\ee

It can be checked  that this family of spaces $\R_{\gamma}$ is a (compact Hausdorff)
projective family, and that  the projective limit of this family is the set of
all, not necessarily continuous, group morphisms from $\R$ to $T$.
This coincides of course, with the dual
group of the discrete group $\R$, which is one of the known characterizations of
the Bohr compactification of the line. Let us denote this space by $\br\equiv{\rm Hom}[\R,T]$.
Being a (commutative) group, $\br$ is naturally equipped with the Haar measure.
 From the above discussion, and in particular from Theorem \ref{teomeba2},
it follows that 
the Haar measure is fully determined by the family of measures
obtained by push-forward, with respect to the  projections:
\be
\label{22}
p_{\gamma}:\br\to\R_{\gamma},\ \ \bar x \mapsto \bar x_{|\gamma},
\ee
where $\bar x_{|\gamma}$ denotes the restriction of $\bar x$ to the subgroup $G_{\gamma}$.
Because each $G_{\gamma}$ is freely generated, it~follows that each space $\R_{\gamma}$ is homeomorphic to a $n$-torus $T^n$, where $n$ is the cardinality of the  set 
$\gamma=\{k_1,\ldots,k_n\}$. 
Furthermore, one can check 
that the push-forward with respect to the projections~(\ref{22}) produces 
precisely the Haar measure on the corresponding torus $T^n$, $\forall\gamma$.
Thus, the measure space $\br$ with corresponding Haar measure can be seen as 
the projective limit
of a projective family of finite dimensional tori, each of which equipped 
with the natural Haar measure.

\section{Measures on Linear Spaces}
\label{esplin}

The infinite dimensional real linear space where a measure can be defined 
in the most natural way is the algebraic dual of some linear space.
We will start by showing that, given any real   linear space
 $E$, its algebraic dual $E^a$ is a projective limit.

Let then $E$ be a real linear space and let us denote by $\cal L$ the
set of all finite dimensional linear subspaces  $L \subset E$. The set $\L$ 
is directed when equipped with the
 partial order relation ``$\geq $'':
\be
\label{bm1}
L\geq L' \mbox{ if and only if } L \supset L',
\ee

Let us consider the family $\{L^a\}_{L\in \L}$ of all spaces dual to
 subspaces  $L\in\L$. For each pair $L,L'$ such~that $L'\geq L$ let
$p_{L,L'}:L'^a \rightarrow L^a$ be the linear transformation  such that  each element of
$L'^a$ is  mapped to its restriction to $L$. The transformations
$p_{L,L'}$ are surjective, since any linear functional on $L$ can be extended
to a linear functional on $L'\supset L$, and the following conditions are satisfied:
\be
\label{bm2}
p_{L,L''}=p_{L,L'}\circ p_{L',L''}\ \ \mbox{\rm for }\ L''\geq L'\geq L.
\ee

It follows that $\{L^a,p_{L,L'}\}_{L,L'\in \L}$ is a  projective family of
linear spaces. Let $E_{\infty}$ be the corresponding  projective limit.
It is clear that $E_{\infty}$ is a linear subspace of the direct product of all spaces
 $L^a$, since the projections  $p_{L,L'}$ are linear. Let 
$\phi$ be a generic element of $E^a$ and $\phi_{|L}$ its restriction
to $L$. Given that, for $L'\geq L$, $\phi_{|L}$ coincides with the
restriction of $\phi_{|L'}$ to the  subspace $L$, one gets a 
linear injective map:
\ba
\label{bm4}
\varpi : \!\!\!\! & E^a & \!\!\!\! \rightarrow E_{\infty} \\
\label{bm5}
& \phi & \!\!\!\! \mapsto \bigl(\phi_{|L}\bigr)_{L\in{\cal L}}.
\ea

On the other hand, the consistency conditions that define $E_{\infty}$
ensure that any element of $E_{\infty}$ defines a linear  functional on
$E$. So, the map $\varpi$ is also surjective, and therefore establishes an
isomorphism  between the  linear spaces $E^a$ and $E_{\infty}$.

Let us consider the measurable projective family  $\bigl\{(L^a,\B _L),
p_{L,L'}\bigr\}$, where $\B _L$ is the Borel $\sigma $-algebra in $L^a$
(recall that $L^a$ is finite dimensional $\forall L$). Let $(E_{\infty},
\B _{\infty})$ be the  measurable  projective limit of this family and  define:
\be
\label{bm6}
\B _{E^a}:=\varpi ^{-1}\B _{\infty}.
\ee

The measurable spaces $(E_{\infty},\B _{\infty})$ and $(E^a,\B _{E^a})$ 
are therefore isomorphic, and we will make no distinction between them. 
The $\s$-algebra $\B _{E^a}$ is the smallest $\sigma $-algebra such that all the 
real functions:
\be
\label{bm7}
E^a\ni \phi\stackrel{\xi}{\longmapsto} \phi(\xi ),\ \xi \in E
\ee
are measurable, i.e.,
\be
\label{bm7a}
\B _{E^a}=\Bs\biggl(\bigcup_{\xi\in E}\xi^{-1}\B(\R)\biggr),
\ee
where $\xi^{-1}\B(\R)$ denotes the family  of inverse images 
of Borel sets of $\R$ by the map (\ref{bm7}).

The fundamental result concerning the  existence of measures on infinite dimensional real linear
spaces is the following \cite{Ya}.
\bt
\label{theobm1}
Any self-consistent family of finite  Borel measures $\mu _L$ on the
subspaces $L^a\subset E^a$ defines a  {($\sigma $-additive)} finite measure on
$(E^a,\B_{E^a})$.
\et
The above result can be presented in a different way,  invoking 
 Bochner's classical theorem.
\bd
\label{defbm1}
Let $E$ be a   real linear space and $\mu $  a finite measure on 
$(E^a,\B_{E^a})$ (if $E$ is finite dimensional, then $E\cong E^a\cong 
\R^n$, $\B_{E^a}$ is the  Borel $\sigma $-algebra in $\R^n$ and $\mu$ is a 
Borel measure).  The Fourier transform, or characteristic function, of the measure  is the  
 (in general complex) function on $E$ 
given by:
\be
\label{bm8}
E\ni \xi \longmapsto \int_{E^a} e^{i\phi(\xi )}d\mu (\phi).
\ee
\ed
\bd
\label{defbm2}
A complex function $\chi $ on a  real  linear space $E$ is said to be of the positive type
if $\sum _{k,l=1}^m c_k \bar c_l\chi (\xi _k-\xi _l) \geq 0,\ \forall
m \in \N ,\ c_1, \ldots, c_m \in \C $ and $\xi_1, \ldots, \xi_m \in E$.
\ed
\bt[Bochner]
\label{theobm2}
A complex function $\chi $ on $\R^n$ is the  Fourier transform  of a finite Borel measure 
on $\R^n$ if and only if it is continuous and of positive type.
The measure is normalized if and only if $\chi(0)=1$.
\et
Bochner's theorem is generalizable to the infinite dimensional situation as follows. 
\bt
\label{theobm3}
Let $\chi $ be a complex function on an 
infinite dimensional real linear space $E$. The function $\chi$ is the  Fourier transform 
of a 
finite 
{measure} on $(E^a,\B_{E^a})$ if and only if it is of the positive type 
and continuous on every  finite dimensional subspace. 
The measure is normalized if and only if $\chi(0)=1$.\et
This result can be proved using  Theorem \ref{theobm1} and
Bochner's theorem. We present next the  essential  arguments. The fact that  the Fourier transform of a measure 
$\m$ on $(E^a,\B_{E^a})$ is necessarily of the positive type on $E$ is a consequence of:
\be
\label{bm9}
\int_{E^a}\Bigl|\sum_k^m c_k e^{i\phi(\xi_k)}\Bigr|^2 d\m(\phi)\geq 0.
\ee

From (\ref{mudvar}) and (\ref{lim12}) one can see that the restriction  of the
Fourier transform of $\m$ to   a finite dimensional subspace $L\subset E$ 
coincides with:
\be
\label{bm10}
\int_{L^a}e^{i\phi_L(\xi)}d\m_L(\phi_L),\ \xi\in L,
\ee
and it is therefore  the Fourier transform of  $\m_L$, hence
continuous. Conversely, a function  $\chi$ of the positive type on $E$ defines, by restriction, 
a family $\{\chi_L\}$ of  positive type functions on the subspaces  $L$:
\be
\label{bm10a}
\chi_L:=\chi_{|L},\ \forall L.
\ee

If $\chi$ is continuous in each $L$ one then have well-defined measures on $L^a$, whose
self-consistency is ensured by (\ref{bm10a}).

\smallskip
To conclude this section, note that for the existence of a measure on $E^a$,  Theorem 
\ref{theobm3} requires only continuity of the characteristic function
on the finite dimensional subspaces. 
Analogously to the situation in finite dimensions, one can expect that a smoother
Fourier transform will produce a measure  supported in proper  subspaces 
of $E^a$. The support of the measure is indeed related to
continuity properties of the   Fourier transform
\cite{Ya,GV,S,GJ,Ri}.
As an  extreme example of this relation, consider the weakest possible topology 
in $E$,  
having only the empty set and 
$E$ itself as open sets. The only continuous functions in this topology are the constant 
functions, and it should be clear that a measure with  constant  Fourier transform is
a Dirac-like measure, supported on the null  element of $E^a$.
In the next section we will discuss two important cases where the
characteristic functions are continuous with respect to a weaker topology   than
the one  defined by continuity in finite dimensional subspaces.
In these cases, the measure is supported in a proper (infinite dimensional) subspace of $E^a$, 
which is equivalent to give a measure on that subspace,
by Proposition \ref{suportegeral} of Section \ref{geral}.

\section{Minlos' Theorem}
\label{B-M-I}

In this  section we consider the relation  between continuity of the  
characteristic function and the support of the corresponding measure, for two situations of interest.

In the first  case the characteristic function is  continuous in a nuclear  
 topology.
In the second case the characteristic function is  continuous with respect to  fixed inner
product.

Let us start by recalling that any family of norms
$\{\|\ \|_{\alpha }\}_{\alpha \in \Gamma }$ in a linear space $E$ 
defines a  locally convex topology
(see, e.g., \cite{Ru}, where the more general case of semi-norms is also considered). 
In fact, one can  take as basis of the topology 
the finite intersections of sets of the form: 
\be
\label{espnuc1}
V(\alpha ,n)=\bigl\{\xi \in E\ |\ \ \|\,\xi \,\|_{\alpha }<1/n\bigr\},\ 
\alpha \in \Gamma,\ n\in \N .
\ee

Also, any family of   norms in the same space  $E$ is 
partially ordered by the   natural order relation
$\|\ \|_{\alpha '}\geq \|\ \|_{\alpha }$ if and only if
$\|\,\xi \,\|_{\alpha '}\geq \|\,\xi \,\|_{\alpha },\ \forall \xi \in E$.
For typical applications, 
it is sufficient to  consider the case where the topology
is defined by a countable and ordered family of  norms,
i.e., we consider sequences of  norms $\{\|\ \|_k\}_{k\in\Ni}$ such that 
$\|\ \|_k\geq\|\ \|_l$, for $k>l$. (In this case the corresponding topology is actually metrizable,  
see, e.g., \cite{KF}.)

For the current application, we restrict attention further to the situation
where the ordered
sequence of  norms $\{\|\ \|_k\}_{k\in\Ni}$ is associated with a
sequence of inner products  $\{\langle \,,\rangle_k\}_{k\in\Ni}$,
$\|\cdot\|_k=\sqrt{\langle \cdot\,,\cdot\rangle_k}$, $\forall k\in\N$. 
With this set-up, let
 $\H_k$ be the  completion of $E$  with respect to the inner product
 $\langle \,,\rangle _k$.
For $k>l$ we have $\H_k\subset \H_l$, since the topology defined by 
$\|\ \|_k$ is stronger. One can show that the topological linear space $E$ defined 
in this way is  complete if and only if
$E=\bigcap_{k=1}^{\infty}\H_k$ \cite{KF}.
\bd
\label{defHS}
An operator $H$ on a separable Hilbert space $\bigl({\cal H},(\ ,\ )
\bigr)$ is said to be a Hilbert-Schmidt operator if given an (in fact any) orthonormal basis
$\{e_k\}$ we have:
\be
\label{espnuc2}
\sum_{k=1}^{\infty} (He_k,He_k) \ \ <\  \infty   .
\ee
\ed
We next define the notion of  nuclear space, following \cite{Ya}. (Note however 
that \cite{Ya}  considers a more general notion, admitting non-countable families of semi-norms,
associated with degenerate inner~products.)
\bd
\label{defespnuc2}
Let  $E=\bigcap_{k=1}^{\infty}\H_k$ be a complete linear space, with respect to the topology defined by an ordered
sequence of  norms $\{\|\ \|_k\}_{k\in\Ni}$  associated with a
sequence of inner products  $\{\langle \,,\rangle_k\}_{k\in\Ni}$.
The space $E$ is said to be  nuclear if $\forall l$
there is $k>l$ and an   Hilbert-Schmidt operator $H$ on 
$\H_k$ such that $\langle \xi,\eta\rangle _l=
\langle H\xi,H\eta\rangle _k$, $\forall \xi,\eta \in \H_k$.
\ed
The most common  examples of nuclear spaces are the following.

\medskip
\noindent {\bf Example 1:} Consider the space $\S $ of rapidly decreasing real sequences 
 $y=(y_n)_{n\in \Ni }$ such that $\lim _{n\to \infty}n^ky_n=0,\ 
\forall k\in \N $, with inner  products:
\be
\label{espnuc3}
\langle y,z\rangle _k=\sum _{n=1}^{\infty}n^{2k}y_nz_n,\ k\in \N .
\ee

For any $k$, the operator $H$ on $\H _k$ (the completion of $\S $ by means of
$\langle \,,\rangle _k$) defined by $(Hy)_n=y_n/n$ is obviously Hilbert-Schmidt.
On the other hand, it is  clear that $\langle\xi,\eta\rangle_k=\langle H\xi,H\eta
\rangle_{k+1}$, and it follows that $\S $ is nuclear.

\medskip
\noindent {\bf Example 2:} The real Schwartz space $\S (\R ^d)$ of
$\cinf$-functions $f$ on $\R ^d$ such that:
\be
\label{espnuc5}
\sup _{x\in\Ri^d}\,\biggl|x_1^{k_1}\ldots x_d^{k_d}
\frac{\partial ^{j_1}}{\partial x_1^{j_1}}\ldots  
\frac{\partial ^{j_d}}{\partial x_d^{j_d}}f(x)
\biggr|<\infty,\ \forall k_1,\ldots,k_d,j_1,\ldots,j_d \in \N 
\ee
is a  nuclear space for an appropriate  sequence of inner
products, whose   topology coincides with the topology
defined by the system of norms  (\ref{espnuc5}) 
 \cite{Ya,GJ} {(see also \cite{BecSen} for more information on the Schwartz~space)}.
\medskip

We present next  the  classical  Bochner-Minlos theorem (whose proof can be found, 
e.g., in \cite{Ya}),
which partially justifies the  relevance of nuclear spaces  in measure theory. 
According to this  result, a characteristic function  which is continuous
in a nuclear space $E$ is equivalent to a measure on the topological dual of $E$. 
 Note that a  linear functional $\f$ on a space of the type $E=\cap_{k=1}^{\infty}\H_k$ is
continuous if and only if it is continuous with respect to (any) one of the inner products
${\langle,\rangle}_k$ \cite{KF}. Equivalently, $\f$ belongs to the
 topological dual $E'$ if and only if $\exists k$ such that $\f\in\H_{-k}$,
where $\H_{-k}$ denotes the (Hilbert space) dual of $\H_k$. So, the
topological dual of a space $E=\cap_{k=1}^{\infty}\H_k$ is a union of
Hilbert spaces, $E'=\cup_{k=1}^{\infty}\H_{-k}$, where
$\H_{-l}\subset\H_{-k}$, for $k>l$. In the  case of the space $\S$ of example 1,
the dual $\S'$ can be seen as the linear space of real sequences  $x=(x_n)$ 
for which there exists $k\in\N$ such that:
\be
\label{dualS}
\sum_{n=1}^{\infty}n^{-2k}x_n^2<\infty .
\ee
In the case of the space $\S(\R^d)$ of example 2, the dual is the space
$\S'(\R^d)$ of tempered distributions, which includes
``Dirac delta functions'' and  derivatives thereof (see, e.g., \cite{GJ}).
\bt[Bochner-Minlos]
\label{theoII2}
Let $E$ be a real nuclear space and $\mu $ a measure on $(E^a,\B_{E^a})$. 
If the characteristic function of the measure is continuous in the nuclear topology,
then the measure is supported on the topological dual  $E'\subset E^a$. 
So, a function of the positive type and continuous on a
nuclear space $E$ defines a 
{measure} on
$(E',\B_{E'})$, where $\B_{E'}:=\B_{E^a}\cap E'$ is the smallest $\s$-algebra
such that all  functions on $E'$ of the type $\phi\mapsto\phi(\xi)$,
$\xi\in E$, are measurable.
\end{theo}
Measure theory in $\S'(\R^d)$ plays a distinguished role in 
applications. The following result  establishes the relation between 
the $\s$-algebra $\B_{\S'(\Ri^d)}$ and the strong topology in $\S'(\R^d)$
\cite{CL} (see~also \cite{Ri}). Recall that the strong  topology
in $\S'(\R^d)$ is generated by the family of semi-norms $\{\rho_A\,|\,
A\subset\S(\R^d)\mbox{\ and bounded}\}$, with $\rho_A(\phi)=\sup_{\xi\in A}
|\phi(\xi)|,\ \phi\in\S'(\R^d)$.
\begin{lem}
\label{SBorel}
The $\s$-algebra $\B_{\S'(\Ri^d)}$ generated by the functions  $\phi\mapsto\phi(\xi)$,
$\phi\in\S'(\R^d)$, $\xi\in\S(\R^d)$, coincides with the  Borel $\s$-algebra 
associated with the strong topology in $\S'(\R^d)$.
\end{lem}
\bc
\label{SBM}
A continuous function of the positive type on $\S(\R^d)$ is equivalent to
a  Borel measure on $\S'(\R^d)$.
\ec
In particular situations,  namely  for Gaussian measures, the  characteristic function
is continuous in a 
topology defined by a single inner  product. In that case 
Minlos' theorem applies.  Minlos'  theorem is presented in the literature
in several different ways (see, e.g.,
\cite{Ya,GJ,Ri,RR}), being most commonly formulated for the case of nuclear spaces.
We start by presenting a more general version, following 
\cite{Ya}, considering next the nuclear space case.
\bt[Minlos]
\label{theo2a}
Let $E$ be a real linear space, $E_1\subset E$ a subspace, $(\,,)$
a inner product in $E$ and $(\,,)_1$ a inner product in $E_1$ such that the
corresponding topology in $E_1$ is stronger than the one induced from 
$\bigl(E,(\,,)\bigr)$. Let $\overline{\bigl(E_1,(\,,)_1\bigr)}$ be the completion
of $E_1$ with respect to $(\,,)_1$. Let  
 $H$ be a Hilbert-Schmidt operator on $\overline{\bigl(E_1,(\,,)_1\bigr)}$ such that:
\be
\label{espnuc6}
(\xi,\eta)=(H\xi,H\eta)_1,\ \ \forall \xi,\eta\in
\overline{\bigl(E_1,(\,,)_1\bigr)}.
\ee
Then, a characteristic function $\chi$ on $E$ continuous  with respect  to
 $(\,,)$, defines a measure supported on the subspace of $E^a$ of those functionals whose
restriction to $E_1$ is continuous with respect to $(\,,)_1$.
\et
In the case of a nuclear space $E$, let us suppose that the 
characteristic function $\chi$ is continuous  with respect  to one of the
inner products  $\langle\,,\rangle_{k_0}$ of the family $\{\langle\,,
\rangle_k\}_{k\in\Ni}$ that defines the topology of $E$. By the very definition 
of  nuclear space, there exist 
$k_1>k_0$ and  a Hilbert-Schmidt operator $H$ such that $\langle\cdot\,,\cdot
\rangle_{k_0}=\langle H\cdot,H\cdot\rangle_{k_1}$. The measure is therefore supported
in $\H_{-k_1}$, the dual of the completion $\H_{k_1}$ of $E$ with respect to
 $\langle\,,\rangle_{k_1}$. More generally we have the following
\bc
\label{1corMinlos}
Let $E$ be a real nuclear space and $(\,,)$, $(\,,)_1$ two inner products in
 $E$ such that the corresponding  topologies are weaker than the 
 nuclear topology. Assume that the   $(\,,)$-topology is weaker than
the  $(\,,)_1$-topology. Let $\overline{\bigl(
E,(\,,)_1\bigr)}$ be the completion of $E$ with respect to $(\,,)_1$ and
  $H$ a Hilbert-Schmidt operator on $\overline{\bigl(E,(\,,)_1
\bigr)}$ such~that:
\be
\label{eq1corMinlos}
(\xi,\eta)=(H\xi,H\eta)_1,\ \ \forall \xi,\eta\in
\overline{\bigl(E,(\,,)_1\bigr)}.
\ee
Then, a characteristic function $\chi$ on $E$ continuous  with respect to
 $(\,,)$, defines a measure supported on the subspace of $E'$ of the functionals 
which
are continuous with respect to $(\,,)_1$.
\ec

\section{Gaussian Measures} 
\label{mgauss}

 In this section we consider Gaussian measures on infinite dimensional real linear spaces, following~\cite{Ya,GV,S}. {In this approach, and following the lines of Section \ref{esplin},
we start with a characteristic function---determined in this case by
an inner product---in a linear space $E$, thus defining the measure initially on the algebraic dual $E^a$. As already mentioned in the Introduction, in other approaches~\cite{Gr,K,Boga},
Gaussian measures are defined directly on topological vector spaces.
The two perspectives are nevertheless equivalent: the algebraic dual $E^a$ is simply
the ``universal home'' for Gaussian measures associated with inner products defined
in $E$. The space where the measure is actually supported is at the end determined by the inner
product itself, regardless of what space one initially considers the measure to be defined in.}

As in  finite dimensions,
Gaussian measures are associated with inner products,  defining the measure's
covariance. (Note that positive semi-definite bilinear forms also give rise
to measures, with the peculiarity that the measure degenerates into a Dirac measure along the null directions. We shall not consider that generalization.)

The fact that the   Fourier transform of a  Gaussian function 
(centered at zero) is also Gaussian allows one to  define 
Gaussian measures on $\R^n$ as follows.
\bd
\label{defmgauss2}
Let $C=(C_{ij})$ be a $n\times n$  positive definite symmetric matrix.
The Gaussian measure $\mu_C$ on $\R ^n$ of covariance $C$ is the
Borel measure whose  Fourier transform is:
\be
\label{mgauss4}
\chi _C(y_1,\ldots ,y_n)=\exp\bigl(-{\textstyle \frac{1}{2}\sum _{i,j}}
C_{ij}y_iy_j\bigr).
\ee
\ed
Using the Lebesgue measure $d^nx$, the 
Gaussian measure of  covariance $C$ is given by:
\be
\label{mgauss5}
d\mu _C(x_1,\ldots ,x_n)=(2\pi )^{-n/2}(\det C)^{-1/2}\exp \bigl(-
{\textstyle \frac{1}{2}\sum _{i,j}}C_{ij}^{-1}x_ix_j\bigr)d^nx.
\ee

A positive definite symmetric matrix is equivalent  to an inner product,
and therefore Gaussian measures on $\R^n$ are determined by
inner  products. One can define
Gaussian measures on infinite dimensional spaces in exactly the same way.
\bd
\label{defmgauss3}
Let $E$  be an infinite dimensional  real linear space and $(\,,)$ an inner
product in  $E$. The  measure on $(E^a,\B_{E^a})$ with  Fourier transform 
$\chi (\xi )=e^{-(\xi ,\,\xi )/2},\ \xi \in E$, is called a
Gaussian measure, of covariance $(\,,)$.
\ed
The existence and uniqueness of the measure are ensured by Theorem
\ref{theobm3} of Section \ref{esplin}. The following
characterization of   Gaussian measures, sometimes taken as definition, is crucial.
\bt
\label{propmgauss1}
A measure $\mu$ on $(E^a,\B_{E^a})$ is Gaussian  if and only if the
push-forward $\mu _{\xi }$ of $\mu $ by the map:
\be
\label{mgauss6}
E^a\ni \phi \longmapsto \phi (\xi )\in \R
\ee
is a Gaussian measure on $\R ,\ \forall \xi \in E$.
\et
The theorem is easily proved. Note first that for any  
Gaussian measure $\mu$ on $E^a$, the push-forward $\mu _{\xi }$ is a
Gaussian measure on $\R$ of covariance
$(\xi,\xi )$. Conversely, let $\mu $ be a measure on $E^a$ such that
$\mu _{\xi }$ is a Gaussian measure on $\R$, $\forall \xi $. Let $c_{\xi }$
be the  covariance of $\mu _{\xi }$. The 
Fourier transform $\chi $ of the measure $\mu $ is then:
\be
\label{mgauss7}
\chi (\xi )=e^{-c_{\xi }/2},\ \forall \xi ,
\ee
where:
\be
\label{mgauss8}
c_{\xi }=\int _{\Ri}x^2d\mu _{\xi }(x)=\int _{E^a}\bigl(\phi (\xi )\bigr)^2d\mu 
(\phi ) .
\ee

On the other hand, it is clear that:
\be
\label{mgauss11}
(\xi_1 ,\xi_2):=\int _{E^a}\phi (\xi_1)\phi (\xi_2)d\mu (\phi )\, ,\ \
\xi_1,\xi_2\in E
\ee
defines an inner product, thus proving that $\chi (\xi )$
is of the required form $\chi (\xi )=\exp \bigl(-(\xi ,\xi )/2\bigr)$.

Expression (\ref{mgauss11}) 
 for the moments of the Gaussian measure of covariance  $(\,,)$ is easily generalized. 
The result is the well-known   Wick's theorem (see \cite{S}). If 
$\xi _1,\ldots,\xi _{2n+1}$ is an odd set of elements of~$E$~
then:
\be
\label{mgauss12}
\int _{E^a}\phi (\xi _1)\ldots \phi (\xi _{2n+1})d\mu =0 .
\ee

If on the other hand  $\xi _1,\ldots,\xi _{2n}$ is an even set of elements of
 $E$ then:
\be
\label{mgauss13}
\int _{E^a}\phi (\xi _1)\ldots \phi (\xi _{2n})d\mu =
\sum _{\mbox  {\rm pairs}}(\xi _{i_1},\xi _{j_1})\ldots 
(\xi _{i_n},\xi _{j_n}) ,
\ee
where $\sum _{\mbox  {\rm pairs}}$ stands  for the sum over all possible ways of 
pairing the $2n$ labels $1,\ldots ,2n$ into $n$ pairs.


Let us note the following. Independently of the linear space $E$ where the  
covariance $(\,,)$ is originally defined, 
a characteristic function of the type  $\chi (\xi )=e^{-(\xi ,\xi )/2}$ 
is always obviously extendable to the Hilbert space completion $\H$ of $E$. 
So, the inner product $(\,,)$, taken as a covariance in the Hilbert space  $\H$,
defines a
Gaussian measure on $\bigl(\H^a,\B_{\H^a}\bigr)$, where $\H^a$ is the algebraic dual
of $\H$ and:
\be
\label{mgauss16a}
\B_{\H^a}=\Bs\biggl(\bigcup_{\xi\in\H}\,\xi^{-1}\B(\R)\biggr) .
\ee

One can show that the  
natural map from $\H^a$ to $E^a$ (defined by the  restriction to $E$ of the  elements
of $\H^a$) is an isomorphism  of measure spaces. 
(From Proposition \ref{mudancavariavel}, the push-forward of the measure
on $\H^a$ is the Gaussian measure on $E^a$ of covariance $(\,,)$, and it follows
that $\H^a\subset E^a$ is a  support of the  Gaussian measure
on $E^a$. To be precise, this map is not strictly measurable, but it establishes 
an isomorphism  between the families of measurable sets 
modulo  zero measure sets, which maps the measure on
 $\H^a$ to the measure on $E^a$.) 
Thus, whenever necessary, one can always assume that the 
covariance  of a  Gaussian measure is defined in a  Hilbert space.

\medskip

\noindent {\bf Example 1:} As in the example of Section \ref{espacoproduto},
let us consider the space $\R^{\Ni}$ of real sequences and the measures
$\m_{\rho}$, given by the product of an  infinite sequence of identical
Gaussian measures on $\R$, each of covariance $\rho$.
Let $\R^{\Ni}_{\rm c}\subset\R^{\Ni}$ be the linear space of those sequences that are 
zero after some order, i.e.,
\be
\label{mgau1}
\R^{\Ni}_{\rm c}:=\{(x_n)\ |\ \exists N_x\in\N \ \ {\rm such \ that}
\ \ x_n=0\ \ {\rm for}\ \ n>N_x\}.
\ee

The space $\R^{\Ni}$ is naturally seen as
 the algebraic dual of $\R^{\Ni}_{\rm c}$, with the action:
\be
\label{mgau2}
x(y)=\sum_n x_ny_n,\ \ x\in\R^{\Ni},\ y\in\R^{\Ni}_{\rm c} ,
\ee
and it is clear that the  product $\s$-algebra
in $\R^{\Ni}$ coincides with $\s$-algebra associated with the interpretation of
$\R^{\Ni}$ as a projective limit. The Fourier transform of the measure  $\m_{\rho}$
is easily seen to be:
\be
\label{mgau3}
\chi_{\rho}(y):=\int_{\Ri^{\Ni}} e^{ix(y)}d\mu_{\rho}(x)=e^{-\frac{1}{2}\rho
\sum_n y_n^2},\,\ \forall y\in
\R^{\Ni}_{\rm c} .
\ee

So, the product measure $\m_{\rho}$ coincides with the Gaussian measure
associated with the inner~product:
\be
\label{mgau4}
\langle y',y\rangle_{\rho}:= \rho\sum_{n=1}^{\infty} y'_ny_n,
\ee
which we assume to be defined
on the   real Hilbert space $\ell^2$ of square summable sequences.
Consider now the space $\S$ of rapidly decreasing sequences 
(Example 1, Section \ref{B-M-I}). Like $\R^{\Ni}_{\rm c}$,  $\S$ is
dense in $\ell^2$ with respect to the topology defined by  
$\langle\,, \rangle_{\rho}$ (which is in fact the  natural  $\ell^2$ topology). 
Moreover, the restriction of $\chi_{\rho}$ to $\S$ is 
continuous in the nuclear topology, since the latter is stronger than the topology
induced in $\S$ from the $\ell^2$-norm.
It then follows from  the Bochner-Minlos theorem 
that the measure $\m_{\rho}$ is supported on the
topological dual $\S'$ of the nuclear space $\S$, for any value of $\rho$. 
Furthermore,  Minlos' theorem allows us
to find  proper subspaces of
$\S'$ that still support the measure.  
Let us now describe this application of  Theorem \ref{theo2a}.
Let then $a=(a_n)$ be an element of $\ell^2$ such that $1\geq a_n>0$, $\forall n$ and let
 $(\,,)_a$ be the inner
product in $\R_{\rm c}^{\Ni}$ given by:
 \vspace{-6pt}
\be
\label{mgauss16h}
(y',y)_a=\sum_n\frac{y'_ny_n}{a_n^2}.
\ee

It is clear that the $(\,,)_a$-topology is stronger than the 
$\ell^2$ topology in $\R_{\rm c}^{\Ni}$. Let $H_a$ be the operator on $\R_{\rm c}^{\Ni}$ 
defined by:
\be
\label{mgauss16i}
(H_ay)_n:=a_ny_n.
\ee

The operator $H_a$ is clearly
Hilbert-Schmidt with respect to $(\,,)_a$, and we have  $\langle y',y\rangle_1=(H_ay',H_ay)_a$. Then, using the usual
characterization of continuous functionals on a 
Hilbert space, it follows from Theorem \ref{theo2a} that the measure $\mu_{\rho}$
is supported on the subspace of $\R^{\Ni}$ of  sequences $x$ \mbox{such that}:
\be
\label{II.11.22}
\sum_{n=1}^{\infty}a_n^2x_n^2<\infty.
\ee

 The subspace
defined by (\ref{II.11.22}) is $H_a^{-1}\ell^2$, i.e., the space of sequences
 $x=(x_n)$ of the  form $x_n=z_n/a_n$, with $z=(z_n)\in\ell^2$.
[Since $a\in\ell^2$, one could be tempted to conclude that the measure is
supported on the space $\ell^{\infty}$ of bounded sequences, but that is 
not the case. It is true that the intersection
 $\bigcap_{a\in\ell^2}H_a^{-1}\ell^2$ of all the spaces 
$H_a^{-1}\ell^2$ coincides with $\ell^{\infty}$, but in fact the space 
$\ell^{\infty}$ is contained in a  zero measure set. There is no contradiction
with $\s$-additivity, since the intersection is not countable.]

Let us remark that given any   Gaussian measure $\m$ of covariance $(\,,)$
in a (real,   infinite dimensional and separable) Hilbert space $\H$, it is always
possible to construct an isomorphism (of measure spaces) mapping the given measure to the 
 Gaussian measure on $\R^{\Ni}$ of the example above, with $\rho=1$ \cite{Ya,S}.
This can be understood as follows.  Let $\{\xi_n\}_{n\in\Ni}$
be an orthonormal basis in $\H$  and consider the map $\theta:\H^a\to\R^{\Ni}$
defined by:
\be
\label{isogau}
\theta(\phi)=\bigl(\phi(\xi_n)\bigr).
\ee

Let $\theta_*\,\m$ be the measure on $\R^{\Ni}$ obtained by push-forward
of $\m$. We then  have (see Proposition \ref{mudancavariavel}):
\be
\label{isofouriergau}
\int_{\Ri^{\Ni}}e^{i\sum_ny_nx_n}d(\theta_*\,\m)=
\int_{\H^a} e^{i\phi\Bigl(\sum_ny_n\xi_n\Bigr)}d\m,\ \ 
\forall (y_n)\in\R^{\Ni}_{\rm c}.
\ee

Given that:
\be
\label{mgau5}
\int_{\H^a} e^{i\phi\Bigl(\sum_ny_n\xi_n\Bigr)}d\m=e^{-\frac{1}{2}\sum_ny_n^2},
\ee
it follows that $\theta_*\,\m$ coincides with the Gaussian measure of the above example,
with $\rho=1$.

\medskip 

\noindent {\bf Example 2:} An  important family of Gaussian measures on $\S'(\R^d)$ is
defined by the following family of inner products:
\be
\label{freemed}
\langle f,g\rangle _{m}=\int f(m^2-\Delta)^{-1}g\,d^dx,\,\ \ f,\ g\in S(\R^d) ,
\ee
where $m\in\R^+$ and $\Delta$ is the Laplacian operator. These measures
are relevant e.g.,\ in quantum  theory  and in certain stochastic processes.

\medskip
We conclude this  section with a   variant of Minlos'  theorem tailored for
Gaussian measures, following immediately from Corollary \ref{1corMinlos},  
Section \ref{B-M-I}.
\bc
\label{teomgauss2}
Let $E$ be a real nuclear space and $(\,,)$, $(\,,)_1$ two inner products in
 $E$ such that the corresponding  topologies are weaker than the 
 nuclear topology. Assume that the   $(\,,)$-topology is weaker than
the  $(\,,)_1$-topology.
Let $\overline{\bigl(
E,(\,,)_1\bigr)}$ be the completion of $E$ with respect to $(\,,)_1$ and
  $H$ a Hilbert-Schmidt operator on $\overline{\bigl(E,(\,,)_1
\bigr)}$ such~that:
\be
\label{mgauss16b2}
(\xi,\eta)=(H\xi,H\eta)_1,\ \ \forall \xi,\eta\in
\overline{\bigl(E,(\,,)_1\bigr)}.
\ee
Then the Gaussian measure of covariance $(\,,)$ is 
supported on the subspace of $E'$ of the functionals 
which
are continuous with respect to $(\,,)_1$.
\ec
Another version of this result, closer to the quantum field theory literature 
\cite{RR,Ri}, 
is the following. 
\bc
\label{rbcor2}
Let $E$ be a   real nuclear space, $(\,,)$ a continuous inner product in $E$ and $\H$
the completion of $E$ with respect to $(\,,)$. Let be $H$ be an injective Hilbert-Schmidt operator 
on $\H$ such that $E\subset H\H$ and $H^{-1}:E\to\H$ is
continuous. Denote by $(\,,)_1$ the inner product
in $E$ defined by $(f,g)_1=(H^{-1}f,H^{-1}g)$. Then the Gaussian measure of covariance $(\,,)$ is 
supported on the subspace of $E'$ of the functionals 
which
are continuous with respect to $(\,,)_1$.
\ec
These two versions are related as follows. Let $H$ be a Hilbert-Schmidt  operator on
$\H$, on the  conditions of Corollary \ref{rbcor2}. The image
$H\H$, equipped with the inner product $(\,,)_1$, coincides with the $(\,,)_1$-completion of $E$. 
Since $H:\H\to H\H$ is
unitary and $H:H\H\to H\H$ is well defined, it follows that $H$ is
 Hilbert-Schmidt on $H\H$. On the other hand we have that:
\be
\label{rbranco5}
(f,g)=(Hf,Hg)_1,
\ee
which shows that the $(\,,)$-topology   is weaker than 
the $(\,,)_1$-topology. Finally, since $(\,,)$ is continuous, the
continuity of  $H^{-1}:E\to\H$ implies that $(\,,)_1$ is also
continuous, and therefore both topologies defined by $(\,,)$ and $(\,,)_1$ are
weaker than the  nuclear topology. All
conditions of Corollary \ref{teomgauss2} are thus~satisfied.

\section{Quasi-invariance and   Ergodicity}
\label{srepci}

We present in this section some concepts 
relevant to the study of transformation properties of measures. 
The notions of quasi-invariance and ergodicity are presented,  together with two important results concerning Gaussian measures. We start by reviewing general notions \cite{RS1}, illustrated with 
straightforward examples.
\bd
\label{rdef1}
Let $\mu _1$ and $\mu _2$ be two  measures on the same  measurable space
$\bigl(M,{\cal B}\bigr)\!$. The measure $\mu _1$ is said  to be absolutely
continuous with respect to $\mu _2$, and we write $\mu _1<\mu _2$, if
$\mu _2(B)=0 \Rightarrow \mu _1(B)=0$.
\ed
As an example, consider the  measures $\m_0$ and $\nu_0$ on
$\R$, where $\m_0$ is the Lebesgue measure and $\nu_0$ is the measure supported
on the interval $I=[0,1]$ defined by $\nu_0(B)=\m_0(B\cap I)$, for any
 Borel set $B\subset\R$. It is clear that $\nu_0<\m_0$, whereas it is not true that
$\m_0<\nu_0$. On the other hand we have for instance the measure
$\m_{\rm Cantor}$ defined by the  Cantor function, which is
supported on the Cantor set (see, e.g., \cite{RS1}). The Cantor set has
Lebesgue measure zero, and therefore the measures $\m_0$ and $\m_{\rm Cantor}$ 
are supported on  disjoint sets.
\bd
\label{eidef3}
Two  measures $\mu _1$ and $\mu _2$  on the same  measurable space $(M,{\cal B})$
are said to be mutually singular, and we write $\mu _1\bot \,\mu _2$, if there exists 
a measurable set $B \in {\cal B}$ such that $\mu _1(B)=0$ and $\mu _2(B^c)=0$, where
$B^c$ is the complement of $B$.
\ed
\bt[Radon-Nikodym]
\label{rteo1}
Let $\bigl(M,{\cal B}\,\bigr)$ be a measurable space and $\mu _1$ and 
$\mu _2$ two  {$\sigma$-finite} measures. The measure $\mu _1$ is absolutely
continuous with respect to  $\mu _2$ if and only if there is a real 
non-negative measurable function
 $f=:d\mu _1/d\mu _2$ on $M$ such that $d\mu _1=fd\mu _2$, i.e.,
$\m_1(B)=\int_Bf\,d\m_2$, $\forall B\in\B$.
\et
The function $d\mu _1/d\mu _2$ in the previous  theorem is said to be 
the Radon-Nikodym derivative.
\bd
\label{rdef2}
Two  measures $\mu _1$ and $\mu _2$  on the same  measurable space 
are said to be mutually absolutely
continuous, or equivalent, and we write
$\mu _1\sim \mu _2$, if $\mu _1<\mu _2$ and $\mu _2<\mu _1$, i.e., if
$\mu _2(B)=0$ if and only if $\mu _1(B)=0$.
\ed
The measures $\m_0$ and $\nu_0$ above are not equivalent. The Gaussian measure
$e^{-x^2/2}\frac{dx}{\sqrt{2\pi}}$, for instance, is equivalent to the
Lebesgue measure. 

The next result  establishes   sufficient and necessary conditions for the
equivalence of two Gaussian measures (centered at the null element) 
 \cite{Ya,S}. 
\bt
\label{medidas equiv}
Let $E$ be  a real infinite  dimensional linear space,
$(\, ,\, )$ and  $(\, ,\, )_1$ two inner  products and $\mu$, $\mu_1$
the  corresponding  Gaussian measures. The measures are
equivalent if and only if the inner product $(\, ,\, )_1$ can be written in the form
 $(f,g)_1=(f,Ag)$, where $A$ is a linear operator defined on the 
$(\, ,\, )$-completion of $E$ such~that:
\begin{enumerate}
\item[(1)] A is bounded, positive and with bounded inverse;
\item[(2)] $A-{\bf 1}$ is Hilbert-Schimdt.
\end{enumerate}
\et
\bd
\label{rdef4}
Let $(M,{\cal B},\mu )$ be a measure space, $\varphi :M \to M$ 
a measurable  transformation and $\mu_{\vf}$ the push-forward of $\mu$.
The measure $\mu $ is said to be invariant  under the  action of
$\varphi $, or $\varphi $-invariant, if $\mu _{\varphi }=\mu $.
If $G$ is a  group of  measurable transformations such that $\mu $ is
invariant for each and every element of $G$, we say that $\mu $ is $G$-invariant.
\ed
As an example,  consider the action of $\R$ on itself, by translations:
\be
\label{srepci5}
x\mapsto x+y,\,\ \forall x\in\R,
\ee
where $y\in\R$. Modulo a multiplicative constant, the Lebesgue measure is the
only  ($\s$-finite) measure on $\R$ which is invariant  under the action of translations
(\ref{srepci5}). {This is a particular case of the well-known 
Haar theorem, which establishes the existence and uniqueness
(modulo multiplicative constants) of (regular
Borel) invariant measures on  locally compact  groups.}

The situation is  radically different in the case of infinite dimensional linear spaces. 
The following argument \cite{HSY} shows for instance that there are no (non-trivial) 
 translation invariant  $\s$-finite Borel measures
in infinite dimensional separable Banach spaces.
Let us  suppose then that such an invariant measure exists, and it does not assign an
 infinite measure to all open balls.
It follows that there is an open ball of radius $R$  with
 finite measure. Since the space is  infinite dimensional, one can find an  
 infinite sequence of  disjoint open balls,  of radius $r<R$, all contained in the first ball. Since by hypothesis the measure is invariant under all  translations, all balls of radius $r$
have the same measure. It follows that this measure is necessarily zero, since all the balls are contained  in the same set, which has  finite measure.
Finally, since the space is separable, it can be covered by a countable
set of open balls of radius $r$, all of them with  zero measure. 
It is therefore proved that the whole space has  zero measure, in  contradiction with the
hypothesis. There are, of course, non $\s$-finite  invariant measures, e.g., the
counting measure which assigns measure 1 to each and every point of the space.
There are also $\s$-finite measures on infinite dimensional spaces which are invariant 
under a restricted set of~translations.

Given a  group of measurable transformations  $G$
on a space  $M$, every $G$-invariant measure $\m$  defines
a  unitary representation $U$ of $G$ in $L^2(M,\m)$, by:
\be
\bigl(U(\vf)\psi \bigr)(x)=
\psi (\vf^{-1}x),\ \ \psi\in L^2(M,\m),\ \vf\in G.
\ee

One can still construct  unitary representations of $G$ using measures
that are not strictly invariant, but instead 
satisfy a weaker condition known as 
quasi-invariance.
\bd
\label{rdef5}
Let $(M,{\cal B},\mu )$ be  a measure space, $\varphi :M \to M$ a
 measurable transformation, and let  $\mu_{\vf}$ denote the push-forward of $\mu$ by $\vf$.
The measure $\mu $ is said to be quasi-invariant  under the action of $\varphi $, or
$\varphi $-quasi-invariant, if $\mu _{\varphi}\sim \mu $. If $G$ is a group of
transformations such that $\mu $ is quasi-invariant for all elements of 
$G$ we say that $\mu $ is $G$-quasi-invariant.
\ed
Regarding the group of translations in $\R$ (or $\R^n$), one can show that any two 
 quasi-invariant measures are equivalent, and therefore equivalent to the Lebesgue measure.
More  generally, when considering  continuous transitive actions
of a   locally compact  group $G$ on a space $M$, there is  a unique
 equivalence class of 
 quasi-invariant measures 
\cite{Ki2}.
\bp
\label{rprop1}
Let $G$ be a  group of measurable transformations on $(M,{\cal B}\,)$ and
$\mu $ a $G$-quasi-invariant measure. The following expression  defines
 a unitary representation $U$ of $G$ in $L^2(M,\mu )$: 
\be
\label{uni}
\bigl(U(\vf)\psi \bigr)(x)=
\biggl(\frac{d\mu _{\vf}}{d\mu }\,(x)\biggr)^{1/2}\,\psi (\vf^{-1}x),
\ee
where $\mu_{\vf}$ denotes push-forward of $\mu$ by $\vf \in G$.
\ep
Going back to the  examples above, one can see that the measure
$e^{-x^2/2}\frac{dx}{\sqrt{2\pi}}$ is quasi-invariant under the action
  (\ref{srepci5}), and thus defines a
unitary representation of translations:
\be
\label{II1.5.2}
\bigl(U(y)\psi\bigr)(x)=e^{-y^2/2+yx}\psi (x-y).
\ee

On the contrary, the measure $\nu_0$, supported on the interval $[0,1]$,
is not  quasi-invariant and cannot possibly provide  a
unitary representation.

Concerning the existence  of   translation quasi-invariant measures on
infinite dimensional spaces,  those are not available either,
in most cases of interest. In particular, one can show the following.
In infinite dimensional   locally convex topological linear spaces
there are no (non-trivial) translation quasi-invariant (i.e., quasi-invariant under all translations) $\s$-finite Borel measures 
  (see  \cite{Ya,F,HSY} and references therein).
Typically,  one finds  situations of quasi-invariance under
a subgroup of the group  of all translations
(like in Theorem  \ref{teorepgau1} below).

We review next some concepts and results from  ergodic theory,
following \cite{Ya,Si,BSZ}. Only finite measures are considered.

\bd[Ergodicity]
\label{eidef2}
Let $(M,{\cal B},\mu )$ be a probability space, where the measure $\mu $ is
$G$-quasi-invariant with respect to a group $G$ of measurable  transformations.
The measure is said to be $G$-ergodic if, for
$B \in {\cal B}$, the~condition: 
\be
\label{conjunto quasi invariante}
\m(B\,\triangle\,\vf B)=0,\ \forall \varphi \in G,
\ee
implies $\m(B)=0$ or $\m(B)=1$.
\ed
In favorable cases of continuous actions in certain
topological spaces, $G$-ergodic measures are supported in a single  orbit of
$G$ (see \cite{Ki2}). In  general we have the following \cite{Ya}.
\bt
\label{ergorbita}
Let  $\m$ be a $G$-quasi-invariant  probability measure  in a measurable space
 $(M,\B)$. The measure is $G$-ergodic if and only if for every
$B\in\B$ with $\m(B)>0$, there exists a countable set 
$\{\vf_k\}_{k\in\Ni}$
of elements of $G$ such that $\m\Bigl(\bigcup_{k=1}^{\infty}\vf_kB\Bigr)=1$.
\et
\bt
\label{eiteo1}
Let $\mu _1$ and $\mu _2$ be two $G$-ergodic measures on the same 
measurable space. Then $\mu _1\sim \mu _2$ or $\mu _1\bot \,\mu _2$. If in 
particular $\mu _1$ and $\mu _2$ are $G$-invariant (and normalized) then
$\mu _1=\mu _2$ or $\mu _1\bot \,\mu _2$.
\et
The following result   establishes also necessary and
sufficient   conditions for ergodicity.
\bt
\label{eiteo2}
Let $\mu $ be a  $G$-quasi-invariant probability measure. The measure
is $G$-ergodic if and only if the only $G$-invariant measurable  (real) functions
are constant, i.e., if and only if the  condition:
$$f(x)=f(\varphi x)\ \ (almost\ everywhere)
\ \ \forall\varphi\in G$$ implies: $$f(x)=constant\ \ (almost\ everywhere).$$
\et
This last  result can be proven  with the following arguments
\cite{Ya,RS1,BSZ}. Suppose that $\m$ is $G$-ergodic. 
Given any invariant real function, the inverse image of any 
 Borel set 
satisfies (\ref{conjunto quasi invariante}), and it is therefore proven 
that ergodicity implies that invariant functions are constant almost everywhere.
Conversely, if a set $B$ satisfies
(\ref{conjunto quasi invariante}), then its characteristic function 
$\chi_B$ (equal to
 1 for $x\in B$ and 0 for $x\not\in B$) is invariant, and the second
 condition on the theorem implies  $\m(B)=0$ or $\m(B)=1$.

 For Gaussian measures the following    important theorem holds \cite{Ya}.
(Essentially, point 1
of  Theorem~\ref{teorepgau1} is what is usually known as the Cameron-Martin theorem. The discussion following Theorem \ref{teorepgau1},
\mbox{as well as} the content of Lemma \ref{lemrepgau1} below, provide in fact illustrations of that theorem.)
\bt
\label{teorepgau1}
Let $(\,,)$ be an inner product in a  real linear space $E$ and $\m$ the corresponding
 Gaussian measure on $E^a$. Let $E^*$ be the subspace of 
$E^a$ of those  functionals that are continuous with respect to the topology defined 
by $(\,,)$, and 
$X$ a subspace of $E^a$, considered as a  subgroup  of the group of
translations in  $E^a$. Then:
\begin{enumerate}
\item[(1)] the measure $\m$ is $X$-quasi-invariant if and only if
$X\subset E^*$,
\item[(2)] the measure $\m$ is $X$-ergodic  if and only if $X$ is
dense in $E^*$.
\end{enumerate}
\et
The following  simplified arguments illustrate point  1 of the  theorem.
Consider the 
Gaussian measure on $\R^n$:
\be
\label{inv1}
d\m(x)=\prod_{j=1}^n e^{-x_j^2/2}\frac{dx_j}{\sqrt{2\pi}}
\ee
and its translation with  respect to $y\in\R^n$. The Radon-Nikodym derivative is:
\be
\label{inv2}
\frac{d\m(x-y)}{d\m(x)}=\exp\biggl(\sum_{j=1}^n x_jy_j\biggr)
\exp\biggl(-\frac{1}{2}\sum_{j=1}^n y_j^2\biggr).
\ee

When considering the  limit  $n\to\infty$, which corresponds to a
measure on $\R^{\Ni}$, one can see  that the  derivative vanishes  unless
 $y=(y_j)_{j\in\Ni}$ is an  element of $\ell^2$. Note that the condition
 $y\in\ell^2$ is actually sufficient 
for  equivalence of the measures, since in that case
$\exp\Bigl(\sum_{j=1}^n x_jy_j\Bigr)$ defines an integrable function on the limit
$n\to\infty$, with respect to the measure
 (\ref{inv1}).
When, on the other hand, one considers  translations by more general elements of 
$\R^{\Ni}$,
one obtains two  (quasi-invariant with respect to $\ell^2$) mutually singular measures.

\section{Gaussian Measures on $\S'(\R^d)$}
\label{temp}
To conclude, we consider the particular, but important case of measures   on
the space of distributions $\S'(\R^d)$ (equipped  with the  Borel $\s$-algebra 
associated with the strong topology--see Lemma \ref{SBorel}, Section \ref{B-M-I}).


Given $g\in S(\R^d)$ one can naturally define an  element of $\S'(\R^d)$, by:
\be
\label{srepci2}
h\mapsto
\int gh\,d^dx,\,\ \forall h\in S(\R^d).
\ee

We will continue to denote  that element by $g$, even if  considered 
as an element of  $\S'(\R^d)$. The inclusion of $\S(\R^d)$ in $\S'(\R^d)$
defined by (\ref{srepci2}) induces an  action of $\S(\R^d)$ in $\S'(\R^d)$,
as a subgroup of the group of translations. Explicitly, given
$g\in S(\R^d)$ we get a measurable transformation in $\S'(\R^d)$:
\be
\label{srepci3}
\phi\mapsto\phi+g,\,\ \forall \phi\in\S'(\R^d).
\ee

Let us say in advance that there are quasi-invariant  normalized  Borel measures,
with respect to the action of $\S(\R^d)$ (\ref{srepci3}).
These measures will simply be called $\S$-quasi-invariant measures.

Let then $\m$ be a $\S$-quasi-invariant measure. From
 Proposition \ref{rprop1}, we then have a  unitary representation of 
(the commutative group) 
$\S(\R^d)$ in 
$L^2(\S'(\R^d),\m)$:
\be
\label{cicli5}
\bigl({\cal V}(g)\psi \bigr)(\phi )=\Bigl(\frac{d\m_g}{d\m}(\phi)
\Bigr)^{1/2}\psi(\phi -g)\, ,\ \ g\in\S(\R ^d),
\ee
where $\m_g$ denotes the push-forward of $\m$ with respect to the map (\ref{srepci3}).

On the other hand, as is typically the case in infinite dimensions,
there are no 
 Borel measures on $\S'(\R^d)$ which remain quasi-invariant under the transitive
action of all translations, i.e., with respect to the natural action of
$\S'(\R^d)$ (seen as a group) on itself \cite{GV,Ya}. [Just like in the discussion 
at the end of the previous section, this  immediately leads to the existence
of  non-equivalent
$\S$-quasi-invariant measures. In fact, given a 
 $\S$-quasi-invariant measure $\m$, it is obvious that the push-forward 
$\m_{\phi_0}(\phi)=\m(\phi-\phi_0)$ defined by any $\phi_0\in\S'(\R^d)$ is
also a $\S$-quasi-invariant measure, and there is 
$\phi_0\bigl(\not\in
\S(\R^d)\bigr)$ such that the two measures are not  equivalent.]

The  simplest examples of $\S$-quasi-invariant measures are Gaussian measures, which we
now consider. In order to simplify the discussion, we impose very strong conditions on 
the measures' covariance. 
Let then  $C$ be a linear continuous bijective operator on $\S(\R^d)$, with  continuous inverse. We~say that $C$ is a covariance  operator $C$ if it is 
bounded, self-adjoint and strictly positive in $L^2(\R^d)$ and if $C^{-1}$, considered
as a   densely defined operator on  $L^2(\R^d)$, is (essentially)
self-adjoint and positive. \mbox{It is} then obvious that the bilinear form:
\be
\label{repgau1}
\langle f,g\rangle _C:=\int f\, Cg\,d^dx,\,\ f, g\in S(\R^d) ,
\ee
in $\S(\R^d)\times\S(\R^d)$ is symmetric, positive and non-degenerate, 
thus defining an inner product $\langle\,,\rangle_C$ in the real linear space
$\S(\R^d)$. A covariance operator $C$ therefore defines a
Gaussian measure, which is  supported in  $\S'(\R^d)$, since the 
 $L^2(\R^d)$-continuity of $C$ ensures that the topology defined by the inner product
$\langle\,,\rangle_C$ is weaker than the  nuclear topology. We will say also that $C$ is 
the measure's covariance, with the understanding that we are referring to an inner product 
of the type (\ref{repgau1}).

Using Theorem \ref{teorepgau1}, one can easily check that these measures are 
$\S$-quasi-invariant and
 $\S$-ergodic. In fact, from the required properties  of the operator $C$ one can write:
\be
\label{repgau3}
\int gh\,d^dx =\int (C^{-1}g)(Ch)\,d^dx =
\langle C^{-1}g,h\rangle _C,\ \ \forall g,h\in\S(\R^d), 
\ee
from what follows  that the functionals  on $\S(\R^d)$ defined by
(\ref{srepci2}) are continuous with respect to 
the  $\langle\,,\rangle _C$-topology. Also, the inclusion of $\S(\R^d)$ in the 
dual  $\S'(\R^d)$  is dense with respect to 
the  $\langle\,,\rangle _C$-topology, since $C^{-1}\bigl(\S(\R^d)\bigr)=\S(\R^d)$.
The conditions of Theorem \ref{teorepgau1}
are therefore satisfied.

In the case of  Gaussian measures, the  Radon-Nikodym derivative 
appearing in (\ref{cicli5}) is easily determined, generalizing the correspondent result in 
 finite dimension:
\begin{lem}
\label{lemrepgau1}
Let $C$ be a covariance operator on $\S(\R^d)$ and $\mu$ the corresponding
measure on $\S'(\R^d)$.
Then:
\be
\label{repgau4}
\frac{d\mu (\phi -g)}{d\mu (\phi )}=e^{-\frac{1}{2}\int g\, C^{-1}g\,d^dx}
\,
e^{\phi (C^{-1}g)},\ \ \forall g\in\S(\R^d).
\ee
\end{lem}
We present next a result \cite{CL} applicable to the important situation of measures that remain invariant under $\R^d$-translations. 
This result characterizes the support  of the measure  in terms of the  local 
behavior of typical  distributions.
To formulate it we need to consider the kernel $\cal C$ of a covariance 
$C$,  defined by:
\be
\label{rbranco5a}
\int f\, Cg\,d^dx =:\int d^dx\,d^dx' f(x){\cal C}(x,x')g(x'),\ \ 
\forall f,g\in\S(\R^d).
\ee

In general, the kernel of the covariance is a  distribution on
$\R^d\times\R^d$. The corresponding measure is invariant under
$\R^d$-translations if and only if ${\cal C}(x,x')={\cal C}(x-x')$.
Let us further recall that a signed measure  on a  measurable space $M$ is
a function on the  $\s$-algebra of $M$ of the form:
\be
\label{rbranco5b}
B\mapsto\int_BFd\nu,
\ee
where $\nu$ is a measure on $M$ and $F$ is an integrable function. In
particular, an (Lebesgue) integrable function on an open set
$U\subset\R^d$ defines a signed measure  on $U$. We will also
say that a distribution $\phi\in\S'(\R^d)$ is a signed measure  in $U\subset\R^d$ if 
there exists a measure $\nu$ on $U$ and an integrable function $F$ 
such that
$\phi(f)=\int_U fFd\nu$, for any function $f\in\S(\R^d)$ supported in $U$.
Then \cite{CL}:
\bp
\label{rbprop1}
Let $\m$ be a  Gaussian measure on $S'(\R^d)$, invariant with respect to
$\R^d$-translations and such that the kernel  $\cal C$ of the covariance 
is not a continuous function. {Then the} support of $\m$ is such that for
$\m$-almost every distribution $\phi\in\S'(\R^d)$
there is no (non-empty) open set
$U\subset\R^d$ on which $\phi$ can be seen as 
a signed~measure.
\ep
\noindent {\bf Example 1:} Let us consider the so-called white noise measures, 
defined by a    covariance proportional to the identity operator , $C=\sigma{\bf 1}$, where  
$\sigma\in\R^+$. Since the  covariance is a scalar, these  measures 
are invariant under
$\R^d$-translations, with  covariance kernel 
 ${\cal C}(x)=\s\delta(x)$, where $\delta$ is the evaluation 
distribution  at $x=0$, i.e.~$\delta(f)=f(0)$, 
$\forall f\in\S(\R^d)$. It follows from the previous
proposition that distributions that can be seen as signed measures
on some open set do not contribute to the measure.
One concludes also  immediately, from Theorem \ref{medidas equiv},  that white noise 
measures are not equivalent to each other, for
$\sigma\neq\sigma'$. Furthermore, from  Theorem \ref{teorepgau1} it follows
that the measures are $\S$-ergodic for any $\s$, and one concludes from
Theorem \ref{eiteo1} that the measures are in fact mutually singular, for 
$\sigma\neq\sigma'$.

\medskip
\noindent {\bf Example 2:} Let us consider again the measures of
Example 2, Section \ref{mgauss}, defined by the covariance~operators:
\be
C_m:=(m^2-\Delta)^{-1} ,
\ee
where $m\in\R^+$ and $\Delta$ is the Laplacian  operator. 
The kernel of $C_m$ is easily found to be:
\be
\label{eqclivre1}
{\cal C}_m(x)=\frac{1}{(2\pi)^d}\int d^dp\,\frac{e^{ipx}}{m^2+p^2}.
\ee

The case $d=1$ ($m\not =0$) corresponds to the path integral for the  quantum harmonic 
oscillator.
(The  particular case $d=1$, $m=1$ corresponds  to the  Ornstein-Uhlenbeck measure.)
For $d>1$ we find measures associated with the path integral formulation 
of quantum field theory. 
For $m=0$ we get the well-known  Wiener measure. (The case  $m=0$, $d<3$, requires 
special care, since the  integral (\ref{eqclivre1}) diverges in the  region
$p\approx 0$. An  appropriate modification leads to the so-called conditional Wiener measure.) 
It is well known that these measures are  supported on continuous functions
for $d=1$ and on  distributions for $d\geq 2$ (see, e.g.,
\cite{GJ,Ri}). In  $d=1$ this result comes from the fact that $(m^2+p^2)^{-1}$
is integrable,  with Fourier transform 
(\ref{eqclivre1}) proportional to
$\frac{1}{m}e^{-m|x|}$. In this situation the test functions in $\S(\R)$ can be replaced
by ``delta functions'', and it makes sense to talk about the two point correlation
function ${\cal C}_m(x,x')$, 
which is proportional to $\frac{1}{m}e^{-m|x-x'|}$. 

\medskip
\noindent {\bf Example 3:} In the canonical approach to the quantization of real scalar field theories in $d+1$
dimensions one looks for \reps\ of the Weyl relations:
\be
\label{weyl}
\V(g)\,\U(f)=e^{i\int fg\,d^dx}\,\U(f)\V(g),
\ee
where $f$ and $g$ belong to $\ss$. What is actually meant by this is 
a pair $(\U,\V)$ of (strongly continuous) unitary \reps\
of the  group $\ss$, satisfying (\ref{weyl}). Any  $\S$-quasi-invariant measure
$\mu$ on $\S'(\R^d)$ produces such a representation. 
In fact, one just needs to consider the Hilbert space
$L^2(\S'(\R^d),\m)$, a 
unitary representation ${\cal V}$ like in (\ref{cicli5}) and a second
unitary representation ${\cal U}$ simply defined by:
\be
\bigl(\U(f)\psi\bigr)(\phi)  =  e^{-i\phi(f)}\psi(\phi).
\ee

Note that whereas the unitary representation ${\cal U}$ is obviously well
defined for any measure, the construction of ${\cal V}$ depends critically
on the $\S$-quasi-invariance of the measure.
It is moreover required that the combined action of $\U$ and $\V$ be
irreducible, which can in turn be seen to be equivalent to $S$-ergodicity of the measure.
Any Gaussian measure therefore satisfies all these criteria.
However, contrary to the situation in finite dimensions,
to  produce a physically meaningful quantization of a given field theory,
the measure must satisfy additional conditions, typically in order 
to achieve a proper quantum treatment of the dynamics, and/or symmetries.
For instance,  the canonical formulation of the free quantum scalar field of
mass $m$ (see, e.g., \cite{Vh} for details) is uniquely associated with the Gaussian measure of covariance:
 \be
{\bf C}_m:=\frac{1}{2}(m^2-\Delta)^{-1/2}.
\ee

\newpage













\end{document}